\documentclass[11pt]{amsart}
\usepackage[applemac]{inputenc} 

\usepackage{amsthm,amssymb,amsbsy,amsmath,amsfonts,amssymb,amscd,mathrsfs}
\usepackage[normalem]{ulem}
\usepackage{caption,subcaption}

\long\def\symbolfootnote[#1]#2{\begingroup%
\def\thefootnote{\fnsymbol{footnote}}\footnote[#1]{#2}\endgroup} 

\usepackage{graphicx}
\DeclareGraphicsExtensions{.jpg}
\usepackage{color}
\usepackage{fullpage}
\usepackage[hypertex]{hyperref}
\newtheorem{theorem}{Theorem}
\newtheorem{corollary}[theorem]{Corollary}
\newtheorem{lemma}[theorem]{Lemma}
\newtheorem{proposition}[theorem]{Proposition}

\theoremstyle{definition} 
\newtheorem{definition}[theorem]{Definition}
\newtheorem{example}[theorem]{Example}
\theoremstyle{remark}
\newtheorem{remark}[theorem]{Remark}

\newcommand{\bt}{\begin{theorem}}
\newcommand{\et}{\end{theorem}}
\newcommand{\bl}{\begin{lemma}}
\newcommand{\el}{\end{lemma}}
\newcommand{\bp}{\begin{proposition}}
\newcommand{\ep}{\end{proposition}}
\newcommand{\bc}{\begin{corollary}}
\newcommand{\ec}{\end{corollary}}
\newcommand{\bdeff}{\begin{definition}}
\newcommand{\edeff}{\end{definition}}
\newcommand{\brem}{\begin{remark}}
\newcommand{\erem}{\end{remark}}
\newcommand{\bex}{\begin{example}}
\newcommand{\eex}{\end{example}}

\newcommand{\be}{\begin{equation}}
\newcommand{\ee}{\end{equation}}

\newcommand{\bi}{\begin{itemize}}
\newcommand{\iii}{\item}
\newcommand{\ei}{\end{itemize}}
\newcommand{\bd}{\begin{description}}
\newcommand{\ed}{\end{description}}

\newcommand{\bqn}{\begin{eqnarray}}
\newcommand{\eqn}{\end{eqnarray}}
\newcommand{\eqnn}{\nonumber\end{eqnarray}}

\newcommand{\nn}{\nonumber}
\newcommand{\ba}[1]{\begin{array}{#1}}
\newcommand{\ea}{\end{array}}

\newcommand{\g}{\gamma}
\newcommand{\al}{\alpha}
\newcommand{\eps}{\varepsilon}

\newcommand{\ph}{\varphi}

\newcommand{\R}{\mathbb{R}}
\newcommand{\N}{\mathbb{N}}

\newcommand{\mb}[1]{\mathbb{ #1 }}
\newcommand{\mc}[1]{\mathcal{ #1 }}

\newcommand{\la}{\langle}
\newcommand{\ra}{\rangle}

\newcommand{\tx}[1]{\mathrm{#1}}

\newcommand{\distr}{\mc{D}}

\newcommand{\norm}[1]{\|#1\|}

\newcommand{\ad}{{\rm ad}}

\begin{document}

\title[Some sub-Finsler  nilpotent structures]{
Sub-Finsler structures from the time-optimal control viewpoint for some nilpotent distributions} 

\author[Barilari]{Davide Barilari}
\address[Barilari]{Universit\'e Paris Diderot, IMJ-PRG, UMR CNRS 7586 - UFR de Math\'ematiques.}
\email{davide.barilari@imj-prg.fr}

\author[Boscain]{Ugo Boscain}
\address[Boscain]{CNRS, CMAP, \'Ecole Polytechnique, Palaiseau,   France,
\& Team GECO, INRIA Saclay}
\email{ugo.boscain@polytechnique.edu}

\author[Le Donne]{Enrico Le Donne}

\address[Le Donne]{Department of Mathematics and Statistics, P.O. Box 35,
FI-40014,
University of Jyv\"askyl\"a, Finland}%
\email{ledonne@msri.org}

\author[Sigalotti]{Mario Sigalotti}
\address[Sigalotti]{INRIA Saclay, Team GECO \&
CMAP, \'Ecole Polytechnique, Palaiseau,   France}
\email{mario.sigalotti@inria.fr}
\date{\today}

\thanks{\emph{Acknowledgements}: This work has been supported by the European Research Council, ERC StG
2009 “GeCoMethods”, contract number 239748.}

\begin{abstract} 
In this paper we 
study the sub-Finsler geometry as a  time-optimal control problem.
In particular,  we consider non-smooth and non-strictly convex sub-Finsler structures associated with 
the Heisenberg, Grushin, and Martinet distributions.
Motivated by problems in geometric group theory, we characterize extremal curves, discuss their optimality, and calculate the metric spheres, proving their Euclidean rectifiability.
\end{abstract}

\maketitle 

 \tableofcontents
\section{Introduction}

Sub-Finsler geometry is a natural generalization of Finsler geometry, sub-Riemannian geometry, and hence Riemannian geometry. 
In this paper, we introduce a very general notion of sub-Finsler structures: at each point of a manifold we consider a subspace of the tangent space endowed with a norm. Such a norm is not necessarily supposed to be strictly convex nor smooth, even away from the origin (the typical example is the $\ell^{\infty}$ norm).
 We will only assume
that this norm changes smoothly  
with respect to the point of the manifold, in a suitable sense (see Section \ref{s:deff}). 
Particularly interesting examples are those  norms  that are ``constant'' with respect to the point, 
since these are the structures that appear in geometric group theory and the theory of 
isometrically homogeneous geodesic spaces. 

Lie groups 
equipped with sub-Finsler structures appear in geometric group theory as asymptotic cones of nilpotent finitely generated groups. 
Indeed, in \cite{pansu} Pansu established that, 
if we look at the Cayley graph of a 
finitely generated 
nilpotent group from afar, such a metric graph
 looks like a
 Lie group endowed with a certain left-invariant geodesic metric. 
 Namely, the sequence of metric spaces $\{ ({\rm Cay}_S(\Gamma),\frac{1}{n}\rho_S)\}_{n \in \N}$ of scaled down Cayley graphs of the nilpotent group $\Gamma$ with generating set $S$ and word metric $\rho_S$ converges in the Gromov--Hausdorff topology \cite{Gromov} towards a 
 Lie group that is stratified and nilpotent
 and is equipped
 with a 
 certain explicit 
left-invariant 
sub-Finsler metric.
We remark that such metrics come from structures that are never sub-Riemannian since the norms are characterized by convex hulls of finitely many points.

Another setting where sub-Finsler structures appear is in 
the study of spaces 
that are isometrically homogeneous,
i.e., 
metric spaces on which the group of isometries acts transitively.
Using the theory of locally compact groups and methods from Lipschitz analysis on metric spaces, \cite{Gleason, Bochner_Montgomery46,mz, b1, b2}, under the additional assumptions of being of finite dimension, locally compact, and the distance being intrinsic, 
it has been proven that 
 these spaces 
are    sub-Finsler manifolds of the following type.
Let $G$ be a connected Lie group and $H$   a compact subgroup. 
Let $\Delta\subseteq T(G/H)$ be a $G$-invariant bracket-generating subbundle of the tangent bundle of the manifold $M:=G/H$.
Consider a function  $F:\Delta \to \R$ that is 
$G$-invariant and, for any $p
\in M$, $F$ restricted to the vector space $\Delta_p$ 
is a norm, i.e., 
it is subadditive, absolutely homogeneous and vanishes only at $0$.
 The sub-Finsler distance (also called Carnot--Carath\'eodory distance)  associated with $\Delta$ and $F$ is defined as 
 $$d(p,q)=\inf {\rm length} _F(\gamma),\qquad \forall p,q\in G/H,$$
 where the infimum is taken along all curves $\gamma$  tangent to $\Delta$ joining $p$ to $q$ and for such curves  ${\rm length} _F(\gamma) := \int F(\dot \gamma)$.
 
Sub-Finsler structures also appear in different applications in control theory, as soon as one considers time-optimal driftless control problems where the controls enter linearly and satisfy polytopic constraints. As an example we can mention time-optimal control problems for three level quantum systems \cite{boscain3level}.



Our paper gives a contribution towards the 
understanding of the geometry 
of sub-Finsler spaces.
Some natural problems are the regularity  of spheres and of geodesics.
For instance, it is an open question whether spheres, i.e., boundaries of metric balls,    are 
rectifiable  from a Euclidean viewpoint.  
Regarding geodesics, 
it is not even known if
any pair of points  
can always be connected by a piecewise smooth
 length-minimizing curve.
 If this is the case, 
 one would like to know if the number of such pieces is
 uniformly bounded.
These are fundamental questions coming directly from the asymptotic study of nilpotent finitely generated groups.
Indeed, there are conjectures about asymptotic expansions for the volume
growth  of balls of large radii that are related to  
the rectifiability of spheres and to 
the above-mentioned regularity of geodesics 
 for the 
asymptotic cone, see \cite{Breuillard-LeDonne1}.

The problem of finding  length-minimizing curves in `constant-type' sub-Finsler geometry    can be locally reformulated as a 
minimum-time 
control problem for a system that is linear in the controls.
Locally, on a manifold $M$, one considers
 $k$ vector fields
  $X_1, \ldots, X_k$ 
  defining the subspace of the tangent space and
a symmetric convex body
  $B\subset\R^k$. The set $B$ identifies the Finsler unit ball $\{u_{1}X_{1}+\ldots+u_{k}X_{k}\,|\,u\in B\}$ in the distribution $\tx{span}\{X_{1},\ldots,X_{k}\}$. 
  The  problem of finding sub-Finsler shortest curves between two points $p_1,p_2$ can be rewritten  as the 
problem of 
finding an absolutely continuous curve $ \gamma(\cdot):[0,T]\to M$, 
together with a measurable function $ u(\cdot): [0,T]\to \R^k$, called control,  
that minimizes the   time $T\geq0$ and satisfies
\begin{equation}\label{problem:intro}
\left\{\begin{array}{ll}
\dot{\gamma}(t)=u_{1}(t)X_{1}(\gamma(t))+\cdots +u_{k}(t)X_{k}(\gamma(t)), &   \\
 {u(t)} \in B, & \\ 
\gamma(0)=p_1,\quad \gamma(T)=p_2.& 
\end{array}
\right.
\end{equation}

The purpose of this paper is to review some techniques of optimal control and to apply them to the study of some low-dimensional key examples: the Heisenberg, the Grushin, and the Martinet 
distributions endowed with norms  whose balls are squares.

Each of these distributions is globally defined by the span of 2 vector fields $X_1,X_2$:
\begin{equation} \label{eq:HGM}
\begin{array}{llll}
 \text{Heisenberg:} &M= \R^3 ,\qquad &X_{1}=\partial_{x}-\frac y2 \partial_{z},\qquad &X_{2}=\partial_{y}+\frac x2 \partial_{z},\\ 
\text{Grushin:} &M=\R^2 ,\qquad &X_{1}=\partial_{x} ,\qquad & X_{2}=x \partial_{y} ,\\ 
\text{Martinet:} &M=\R^3 ,\qquad &X_{1}=\partial_{x}+y^{2} \partial_{z},\qquad &X_{2}=\partial_{y}.
\end{array}
\end{equation}
The existence of time-minimizers for the problem \eqref{problem:intro} in these three cases is a classical consequence of Filippov's theorem.

We consider the $\ell^\infty$ and $\ell^1$ norm with respect to the above $X_1,X_2$. More precisely, 
for every $v\in {\rm span}\{X_1(q),X_2(q) \}\subseteq T_qM$, we consider the two quantities
\begin{eqnarray*}
	\norm{v}_\infty &:=& \min \{ \max( |w_1|,|w_2| )   \;:\;   v= w_1X_1(q)+w_2 X_2(q)\} ,   \\
	\norm{v}_1 &:=& \min \{  |w_1|+|w_2|    \;:\;   v= w_1X_1(q)+w_2 X_2(q)\}   .
\end{eqnarray*}
Notice that for an arbitrary structure defined by vector fields $X_1,X_2$ 
the  $\ell^1$ norm with respect to  $X_1,X_2$ is the 
$\ell^\infty$  norm with respect to  $\frac{1}{2} (X_1+X_2)$ and $\frac{1}{2} (X_1-X_2)$. 
In the case of Heisenberg group it is easy to see that the $\ell^1$ structure coincides with the $\ell^\infty$ structure in a new system of coordinates. For the cases of Grushin and Martinet we obtain, up to a multiplicative constant, the vector fields 
\begin{equation} \label{eq:HGM22}
\begin{array}{llll}
\text{Grushin:} &M=\R^2 ,\qquad &X_{1}=\partial_{x}+x \partial_{y} ,\qquad & X_{2}=\partial_{x} -x \partial_{y} ,\\ 
\text{Martinet:} &M=\R^3 ,\qquad &X_{1}=\partial_{x}+\partial_{y}+y^{2} \partial_{z},\qquad &X_{2}=\partial_{x}-\partial_{y}+y^{2} \partial_{z}.
\end{array}
\end{equation}

Regarding our first result, recall that a {\em bang-bang trajectory} is a finite concatenation of  curves, called \emph{arcs}, corresponding to a control that is constant with values in $\{ (1,1), (1, -1), (-1,1), (-1,-1)\}$. 

\begin{theorem} Consider the sub-Finsler structures for the Heisenberg, Grushin and Martinet distribution defined by $\ell^\infty$ norm with respect to vector fields \eqref{eq:HGM} or \eqref{eq:HGM22}.
Then the length-minimizing  trajectories are curves of two types:
\begin{itemize}
\item[(i)] one component of the control is constantly equal to $1$ or $-1$, 
\item[(ii)]  bang-bang trajectory.
\end{itemize}
Moreover, the length-minimizing  trajectories that are not of type (i) have at most 7 arcs.
In addition, all curves of type (i) are length-minimizers.
\end{theorem}
We remark that type (i) and type (ii) are not mutually exclusive. Moreover, it turns out that for every trajectory of type (i) there exists a length-minimizing trajectory of type (ii) connecting the same two points.
As a corollary, we deduce that
any pair of points  
can  be connected by an optimal bang-bang trajectory with at most 7 arcs.

The proof is based on the classical Pontryagin Maximal Principle, for the description of extremal trajectories, i.e., trajectories that satisfy a first-order optimality condition. The bound on the number of arcs for optimal trajectories is obtained via second-order optimality conditions proposed by Agrachev and Gamkrelidze in \cite{Agrachev_Gamkrelidze_1990}.

\medskip
 
 In all the three cases we give a complete description of the sub-Finsler spheres. 
The case of the Heisenberg group was already studied in \cite{Breuillard-LeDonne1} with metric methods and it was proved that the sub-Finsler sphere is Euclidean rectifiable. Here we recover the shape of the sub-Finsler spheres and in addition we obtain the sub-Finsler front, i.e., the set of endpoints of extremals at a fixed time.

We remark that 
apart from the  
Heisenberg distribution, which comes form a left-invariant structure with respect to a group law, the other two examples that we study are not homogeneous structures.
Nonetheless, they come from 
projections of homogeneous structures on groups. Namely, the Grushin plane is a right/left quotient of the Heisenberg group equipped with a right/left-invariant strucure
and 
the Martinet space is a right/left quotient the Engel group  (which is the simplest stratified  group of rank 2 and step 3) equipped with a right/left-invariant strucure.
Consequently, given a length-minimizer in Grushin (resp.~Martinet) the curve in Heisenberg (resp.~Engel) with the same control is length-minimizer as well.
A feature of being quotients of nilpotent groups, with the induced  
projected structure, is that the distribution can be defined by vector fields
that generate a nilpotent Lie algebra.
Nilpotency 
simplifies considerably the problem. For example, in our case bang-bang trajectories  have piecewise polynomial coordinates. 

Thanks to the previous description of length-minimizing trajectories we are able to parametrize the spheres  for the Grushin plane and the Martinet distribution, and we obtain as a consequence the following result.

\bt Consider the sub-Finsler structures for the Heisenberg, Grushin and Martinet distribution defined by $\ell^\infty$ and $\ell^1$ norm with respect to vector fields \eqref{eq:HGM}. Then the sub-Finsler spheres are Euclidean rectifiable, semianalityc and homeomorphic to Euclidean spheres. 
\label{222}\label{223}
\et

The semi-analyticity of spheres 
is interesting since it does not hold 
  for the sub-Riemannian Martinet sphere, as proved in \cite{ABCK97}. 
Indeed in that paper it is proved that  the sub-Riemannian Martinet sphere is  even not sub-analytic.


We mention that there are a few other works that consider the view point of sub-Finsler geometry.
A part from the previously mentioned ones,
in the papers
\cite{clellandmoseley06, clellandmoseley07}
the authors study the sub-Finsler geometry, such as geodesics and 
rigid curves,
in three-dimensional manifolds and in 
 Engel-type manifolds. However, in those papers 
 there is an assumption that is classical in Finsler geometry:
 the norm is assumed to be smooth outside the zero section and strongly convex. 
 The present paper deals mainly with the case where these assumptions are not satisfied. 
Another notable paper is
 \cite{cowlingmartini13}, in which the authors study the 
sub-Finsler  geometry associated with the
solutions of evolution equations given by first-order differential operators, providing one more setting where sub-Finlser geometry appears naturally.
In a paper in preparation  \cite{LeDonne-NicolussiGolo},
the authors study the Euclidean Lipschitz regularity of arbitrary left-invariant distances in a family of homogeneous groups, including the Heisenberg group. 
A significant remark  is that the sphere with respect to the sub-Riemannian distance in 
the product of the Heisenberg group with the real line is a Lipschitz domain, while the  sphere for the $\ell^1$ sub-Finsler structure on the same group admits a cusp, as observed in \cite{Breuillard-LeDonne1}.

\subsection{Definitions} \label{s:deff}

%
 
A function  on $\R^k$, $k\in \N$, is
a {\em norm} if it is subadditive, absolutely homogeneous and vanishes only at $0$.

A {\em    sub-Finsler structure}  (trivialized and of constant-type norm) of rank at most $k$ on a smooth manifold $M$  is 
a pair 
$(  f,  \norm{\cdot} )$ where $\norm{\cdot}$ is a norm on $\R^k$ and $f:  M\times \R^k \to TM$ is
a smooth morphism of bundles  such that $f(\{p\}\times \R^k) \subseteq T_p M$, for all $p\in M$.


With every sub-Finsler structure we associate the distribution $\distr= f(E)$ and a norm on 
$\distr
$ defined by
\begin{equation}\label{norma}
 \norm{v}_{\rm sF}= \inf \{ 
 \norm {w}
 \;:\; f(p,w)=v\},\qquad \text{ for all } v\in \distr_p.
 \end{equation}
 
A  distribution $\distr\subseteq TM$ is 
 \emph{Lie bracket generating}
 if 
 $({\rm Lie} (\Gamma(\distr))) _ p = T_p M$, for all $p\in M$.
 Here $\Gamma(\distr)$ is the collection of smooth vector fields tangent to $\distr$ and, given a family $\mathcal F$ of vector fields, we denote by  $\rm Lie(\mathcal F)$ and $\mathcal F_p$    the Lie algebra generated by $\mathcal F$  and the evaluation of the elements of $\mathcal F$ at a point $p$, respectively.

\begin{remark}

The   norm on $\R^k$, defined as above, 
in not necessarily Finsler in the classical sense, since it is not smooth away from the origin. 
As a consequence, even when 
$\distr= TM$ the function $v\mapsto  \norm{v}_{\rm sF}$ does not necessarily endow $M$ with a  Finsler structure in the classical sense.

%
%
%
%
%

The notion of sub-Finsler structure introduced above, when some non-smoothness of the norm is allowed, 
can be seen as  a particular case of the following more general class:
%
A {\em  partially smooth sub-Finsler structure}   on $M$  is 
a triple $(E, \norm{\cdot}_E, f)$ where    $E$ is a vector bundle over $M$, $\norm{\cdot}_E$ is a partially smooth Finsler structure on $E$ (defined following Matveev and Troyanov \cite{Matveev-Troyanov}), and $f: E \to TM$ is
a smooth morphism of bundles  such that $f(E_p) \subseteq T_p M$, for all $p\in M$.
The norm on the induced distribution can be defined in analogy with  \eqref{norma}, replacing $\|w\|$ by
$\norm{(p,w)}_E .$
\end{remark}

Given a sub-Finsler structure
with distribution $\distr$ and norm $\norm{\cdot}$ we say that an absolutely continuous  curve $\g:[0,T]\to M$ is
 \emph{horizontal} if  $\dot \g (t) \in \distr_{\g(t)}$ and in this case its \emph{lenght} is defined by
$$\ell(\g)=\int_{0}^{T}\norm{\dot \g (t)}_{\rm sF} dt.$$
We can then define the induced distance
$$d(p,q)=\inf\{\ell(\g)\;:\; \g:[0,T]\to M \text{ horizontal and }  \g(0)=p,\g(T)=q\},$$
which is well defined and finite if the distribution $\distr$ is {Lie bracket generating}.

%
%

\section{Sub-Finsler geodesics as minimizers of a time-optimal control problem} 
%


Let $M$ be a smooth manifold and $(f, \norm{\cdot})$ a sub-Finsler structure on $M$.
Notice that the  bundle morphism  $f: M\times \R^k \to TM$    
 determines $k$ vector fields
  $X_1, \ldots, X_k$ defined by $X_i(p) = f (p , e_i)$, where $e_{1},\ldots,e_{k}$ is an orthonormal basis for $\R^{k}$.
(Conversely,
given any $k$ vector fields
  $X_1, \ldots, X_k$ on   
     $M$ 
     there exists a unique bundle morphism $f: M\times \R^k \to TM$
     for which $X_i(p) = f (p , e_i)$.)

 The norm $\norm{\cdot}$ identifies the set
 $$B:= \{ w\in \R^k \;:\;  \norm{w}\leq 1\},$$
which is a closed, convex, centrally symmetric, and with the origin in its interior.
(Conversely, any such a set is the closed unit ball of a norm on $\R^k$.)

The  problem of finding sub-Finsler geodesics, i.e., curves that minimize the length between two points $p$ and $q$, can be reinterpreted as a time-optimal control problem, that is the problem of minimizing the time $T\geq0$ for which
there exist $\gamma:[0,T]\to M$   absolutely continuous and $u: [0,T]\to \R^k$ measurable such that
\begin{equation}\label{problem}
\left\{\begin{array}{ll}
\dot{\gamma}(t)=u_{1}(t)X_{1}(\gamma(t))+\cdots +u_{k}(t)X_{k}(\gamma(t)), & \text{ for almost every } t\in [0,T], \\
 {u(t)} \in B, & \text{ for almost  every } t\in [0,T],\\ 
\gamma(0)=p,\quad \gamma(T)=q.& 
\end{array}
\right.
\end{equation}


Notice that the control function $u$ 
might not be uniquely determined by the trajectory $\gamma$, since the vector fields might not be linearly independent at every point.  
However, given the control $u$ there exists a unique trajectory $\gamma$ satisfying 
$\dot{\gamma}(t)= \sum_{j=1}^k u_{j}(t)X_{j}(\gamma(t)) $ and
$\gamma(0)=p$.

\subsection{Hamiltonian formalism and Pontryagin Maximum Principle}

If a pair $(\gamma(\cdot),u(\cdot))$ is a time-minimizer for \eqref{problem}, then it satisfies the first-order necessary conditions given by the Pontryagin Maximum Principle (PMP). Here we state a suitable version of the PMP for time-optimal control problem on a manifold $M$ (see, for instance, \cite[Corollary 12.12]{Agrachev_Sachkov}).

Define the {\em Hamiltonian}
\bqn
\mc{H}(\lambda,p,u):= \la \lambda, f(p,u)\ra =  \sum_{i=1}^{k} u_{i}\la \lambda, X_{i}(p)\ra ,
\eqn
for $\lambda \in T^*_pM$, $p\in M$, and $u\in \R^k$.
For every $u\in \R^k$, let $\vec{\mc{H}}(\cdot, \cdot, u)$ be the vector field on $T^*M$ uniquely determined by 
the relation
$$ \sigma ( \cdot,  \vec{\mc{H}}(\lambda , p, u))= d_{(\lambda,p)} \mc{H}(\lambda,p,u), $$
where $\sigma$ is the canonical symplectic form on $T^*M$.

Define the {\em maximized Hamiltonian} 
\begin{equation}\label{maximized Hamiltonian}
H(\lambda,p):=\max \{ \mc{H}(\lambda,p,u) \;:\;    {u} \in B\}  .
\end{equation}

\bt[PMP]\label{PMP} Let $(\gamma(\cdot),u(\cdot))$ be a time-minimizer for   Problem \eqref{problem}.
Then there exist an  absolutely continuous function $\lambda:[0,T]\to T^{*}M$  and a constant $\lambda_{0}\geq 0$ such that
\bi
\iii[(i)] 
$\lambda(t)\in T^{*}_{\gamma(t)}M \setminus \{0\},$ for every $t\in [0,T]$,
\iii[(ii)] the pair $(\lambda(t),\gamma(t))$ satisfies the Hamiltonian equation 
$$(  \dot{\lambda}(t)  , \dot{\gamma}(t)) = \vec{\mc{H}}(\lambda(t) , \gamma(t), u(t))),\qquad  \text{ for almost every } t\in [0,T],$$
which, in canonical coordinates, is 
$$\dot{\lambda}(t)=-\frac{\partial \mc{H}}{\partial p}(\lambda(t),\gamma(t),u(t)), \quad \dot{\gamma}(t)=\frac{\partial \mc{H}}{\partial \lambda}(\lambda(t),\gamma(t),u(t)), \qquad \text{ for almost every } t\in [0,T],$$
\iii[(iii)] $\mc{H}(\lambda(t),\gamma(t),u(t))=H(\lambda(t),\gamma(t))=\lambda_{0}$, for almost every  $t\in [0,T]$. 
\ei
\et

If $\lambda(\cdot),\gamma(\cdot)$ satisfy for some $u(\cdot)$ and $ \lambda_{0}$ the conditions (i), (ii), (iii) of  Theorem \ref{PMP},  we say that
 $(\lambda(\cdot),\gamma(\cdot))$ is an {\em extremal pair}, that
 $\gamma(\cdot)$ is an {\em extremal trajectory},
  and that $\lambda(\cdot)$ is an {\em extremal lift} of $\gamma(\cdot)$.


For every vector field $Y$, if   $(\lambda(\cdot),\gamma(\cdot))$ is an extremal pair, then the function
$t \mapsto
 \la \lambda(t) ,Y(\gamma(t))   \ra$
is absolutely continuous and   its derivative satisfies
\begin{equation}\label{derivata-switch}
\dfrac{d}{d t}\la \lambda(t) ,Y(\gamma(t))   \ra  = \la \lambda (t) , \sum_{j=1}^k u_j(t) [X_j, Y](\gamma(t)) \ra,
\end{equation}
for almost every $t$ as it follows from the next classical computation.
In  canonical coordinates,  thanks to point (ii) in Theorem \ref{PMP}, one has
\begin{eqnarray*}
\dfrac{d}{d t}\la \lambda(t) ,Y(\gamma(t))   \ra  &=&
\dfrac{d}{d t}   ( \lambda(t)^T Y(\gamma(t))      \\
&=&\left(-\frac{\partial \mc{H}}{\partial p}(\lambda(t),\gamma(t),u(t))  \right)^T Y  (\gamma(t))  
+\lambda(t)^T   \frac{\partial Y}{\partial p}(\gamma(t))    \frac{\partial \mc{H}}{\partial \lambda}(\lambda(t),\gamma(t),u(t))  
\\
&=&
-  \lambda(t)^T   \frac{\partial f}{\partial p}(\gamma(t), u(t))   Y(\gamma(t))  + \lambda(t)^T   \frac{\partial Y}{\partial p}(\gamma(t))   f(\gamma(t), u(t))
\\&   =&
\la \lambda(t), [f(\cdot , u(t)), Y] (\gamma(t))   \ra .
\end{eqnarray*}
 
\subsection{Second-order optimality conditions}
Our aim is to recall necessary conditions for the optimality of
an extremal   trajectory whose corresponding control is piecewise constant. 
We refer to \cite{Agrachev_Gamkrelidze_1990}. (See also \cite{AgrachevSigalotti,Sigalotti_JDCS}.)

\begin{theorem}\label{thm2nd}
Let $M$ be a smooth manifold and $f: M\times \R^k \to TM$ a 
sub-$\ell^\infty$ structure   on $M$.
Let  $(\gamma(\cdot),u(\cdot))$ be an extremal pair for Problem \eqref{problem}
 and let $\lambda(\cdot)$ be an extremal lift of $\gamma(\cdot)$.
Assume that $\lambda(\cdot)$ is the unique  extremal lift of $\gamma(\cdot)$, up to multiplication by a positive scalar.
 Assume that there exist
 $0=\tau_0<\tau_1<\tau_2<\cdots<\tau_K<\tau_{K+1}=T$ and $u^0,\dots,u^K\in \R^k $ such that   
 $u(\cdot)$ is constantly equal to $u^j$ on $(\tau_{j},\tau_{j+1})$, for $j=0, \ldots, K$.
  
Fix $j=1, \ldots, K$.  For $i=0, \ldots, K$ let $Y_i=f(\cdot , u^i)\in {\rm Vec}(M)$ 
and 
 define recursively the operators
$P_j=P_{j-1}={\rm id}_{{\rm Vec}(M)}$,
$$P_i=   P_{i-1}  \circ   e^{  (\tau_{i}-\tau_{i-1}  ) \ad (Y_{i-1})}  ,\quad\forall  i \in \{j+1, \ldots, K\},$$
$$P_i= P_{i+1}  \circ   e^{ - (\tau_{i+2}-\tau_{i+1}  ) \ad( Y_{i+1})}  , \quad\forall i\in \{ 0,\ldots, j-2\}.$$
Define the   vector fields  
$$Z_{i}=  P_i(Y_i)  
,\quad\forall  i \in \{0, \ldots, K\}.$$
Let $Q$ be the quadratic form
\be\label{2nd} Q(\alpha)=\sum_{0\le i < l\le K}\alpha_i\alpha_l \la \lambda(\tau_j),
[Z_i,Z_l](\gamma(\tau_j))\ra
\,,
\ee
defined on the space
\be\label{space_of_alphas}
W=
\left\{\alpha=(\alpha_0,\dots,\alpha_K)\in \R^{K+1}\Big| 
\sum_{i=0}^K\alpha_i=0,\ \sum_{i=0}^K \alpha_i Z_i(\gamma(\tau_j))=0\right\}.
\ee
If 
$Q$ is not negative semi-definite, i.e., if there exists $\alpha \in W$ such that $Q (\alpha)> 0$, then $\gamma(\cdot)$ is not time-minimizing.
\end{theorem}

\section{Sub-$\ell^\infty$ structures} 
A choice of norm in $\R^k$ that is of particular interest is the
  $\ell^\infty$-norm, that is,
$$\norm {w}=|{w}|_\infty:=\sup_{i=1,\ldots,k} |w_{i}|.$$
When the norm in the definition of  sub-Finsler structure is the  $\ell^\infty$-norm
 we speak about  {\em sub-$\ell^\infty$ structure}. 
The corresponding time-optimal control problem \eqref{problem}
 rewrites as
 \begin{equation}\label{problem2}
\left\{\begin{array}{ll}
\dot{\gamma}(t)=u_{1}(t)X_{1}(\gamma(t))+\cdots +u_{k}(t)X_{k}(\gamma(t)), & \text{ for almost every } t\in [0,T], \\
 {|u_i(t)|} \leq 1, & \text{ for   every }i\in \{ 1,\ldots,k\} , \;t\in [0,T],\\ 
\gamma(0)=p,\quad \gamma(T)=q.& 
\end{array}
\right.
\end{equation}
 The maximized Hamiltonian \eqref{maximized Hamiltonian} is 
\begin{equation}\label{maximized Hamiltonian2}
H(\lambda,p)=|\la \lambda, X_{1}(p)\ra|+\cdots +|\la \lambda, X_{k}(p)\ra|.
\end{equation}

\subsection{Switching functions, singular, abnormal, and regular arcs}

With every extremal pair $(\lambda(\cdot),\gamma(\cdot))$ and every $j=1,\dots, k$
we associate 
   the \emph{switching function} 
   $$t \mapsto\varphi_{j}(t):=
   \la \lambda(t) ,X_{j}(\gamma(t))\ra
   .$$
 By formula
 \eqref{derivata-switch}
 we have that
 \begin{equation}\label{derivata-switch2}
\dot \varphi_j (t) = \la \lambda (t) , \sum_{i=1}^k u_j(t) [X_i, X_j](\gamma(t)) \ra\qquad\mbox{for almost every $t$.}
\end{equation}

The maximality condition (iii) of the PMP and 
 \eqref{maximized Hamiltonian2}
 imply  that
\begin{equation}\label{lambda_zero}
|\varphi_{1}(t)|+\cdots +|\varphi_{k}(t)|=\lambda_{0}, \quad \text{ for all } t
\end{equation}
and that, 
for all $j=1,\ldots, k$ and almost every $t$,
\begin{equation}\label{segno}  \varphi_{j}(t)\neq 0 \implies u_{j}(t)= {\rm sign}\,\varphi_{j}(t) .
\end{equation}

The restriction of an extremal pair $(\lambda(\cdot),\gamma(\cdot))$  to some open nonempty interval $I\subset [0,T]$ is called
\bi
\iii[(i)]  an \emph{abnormal arc} if $\varphi_{j}(t)\equiv 0$ on $I$ for all $j=1,\ldots, k$;
\iii[(ii)] a \emph{$\varphi_{j}$-singular arc} if $\varphi_{j}(t)\equiv 0$ on $I$;
\iii[(iii)]   a \emph{regular arc}  if $\varphi_{j}(t)\neq 0$ for every $t\in I$ and for every $j=1,\ldots, k
$;
\iii[(iv)] a \emph{bang arc} if the control $u(\cdot)$ associated with the trajectory 
is constant and takes values in $\{ 1, -1\}^k$.
\ei


Notice that a regular arc is a bang arc, but the converse is not true. Indeed, bang arcs can be singular (see Section \ref{s:martinet}).

A {\em bang-bang trajectory} is a curve corresponding to a  control that is piecewise constant with values $\{ 1, -1\}^k$. 
In particular, a concatenation of regular arcs is a bang-bang trajectory, called \emph{regular bang-bang trajectory}. When not specified otherwise arcs are assumed to be maximal, meaning that the restriction of the extremal pair to strictly larger open intervals is not an arc.

\begin{remark}\label{rem:abnormal}
An arc is abnormal   if and only if it is $\varphi_{j}$-singular   for all $j=1,\ldots, k$ and if and only if
$\lambda_0= 0$.
The latter equivalence follows from \eqref{lambda_zero}.
In particular,
if a trajectory contains an abnormal arc then  the whole trajectory is an abnormal arc.
\end{remark}

\section{Heisenberg group}

In this section we provide a description of the time-minimizing trajectories in the sub-$\ell^\infty$   Heisenberg group. These results have been previously obtained 
in \cite{Breuillard-LeDonne1}
using methods of metric geometry.
The aim of this section is to illustrate how to exploit  the geometric-control tools presented in the previous sections to recover such results.

We consider the  sub-$\ell^\infty$ structure on the Heisenberg group $\mb{H}\simeq \R^3$
determined by the vector fields 
\bqn\label{eq:campiH}
X_{1}=\partial_{x}-\frac y2 \partial_{z},\qquad X_{2}=\partial_{y}+\frac x2 \partial_{z}.
\eqn
Let us introduce the vector field $X_3=\partial_{z}$, which satisfies $[X_1,X_2]=X_3$ and  $[X_1,X_3]=[X_2,X_3]=0$.

We use the notation from the previous section.
Formula \eqref{derivata-switch2} gives immediately
\bqn \label{eq:pmpphi}
\dot{\varphi}_{1}=-u_{2}\varphi_{3},\qquad \dot{\varphi}_{2}=u_{1}\varphi_{3}, \qquad \dot{\varphi}_{3}=0,
\eqn where  $\varphi_3(t):=\la \lambda(t) ,X_3(\gamma(t))\ra$. 
 Notice that, since $X_1, X_2, X_3$ are linearly independent at every point,  this is a reformulation in coordinates of the vertical part of the Hamiltonian system of the PMP.

We 
characterize here below the abnormal, singular, and regular arcs for the associated time-optimal control problem.
First, we show that there is no nontrivial abnormal trajectory. Next, we describe   the structure of regular and singular arcs,
 showing 
that every nonconstant extremal trajectory   is 
either a singular arc  
or  a concatenation of regular arcs. 
Finally, we give a bound on the maximal number of regular arcs of a time-minimizer.

\subsection{Abnormal arcs}
\bl \label{lemma-abnormal-H}
 The only abnormal arcs on $\mb{H}$ are the constant curves.
 Consequently, no minimizer joining two distinct points is abnormal.
\el
\begin{proof}  

From  Remark \ref{rem:abnormal}, we have
  $\varphi_{1}(t)=\varphi_{2}(t)=0$ for all $t$. 
By non-triviality of the covector $\lambda(\cdot)$, we deduce that $\varphi_3(t)\ne 0$ for every $t$. 

By the first two equations in  \eqref{eq:pmpphi}, we  get $u_{1}(t)=u_{2}(t)=0$ for almost every $t$. 
Then the trajectory is constant and does not minimize the  time. 
\end{proof}

\subsection{Singular arcs}

\bl \label{lemma-singular-H}
 On $\mb{H}$
the nonconstant trajectories that have   singular arcs are 
exactly those
for which
there exists $j\in \{1,2\}$ such that
$u_j$ is constantly equal to $1$ or $-1$.
All of them consist of a single singular arc and are time-minimizers.
\el
\begin{proof}
In what follows the roles of $u_1$ and $u_2$ are interchangeable. 
Consider an extremal trajectory that  is not trivial and is $\varphi_1$-singular when restricted to an interval $I$, i.e.,
 $\varphi_{1}\equiv 0$ on   $I$.
Because of Lemma~\ref{lemma-abnormal-H}, the trajectory does not have abnormal arcs, i.e., $\lambda_0\neq 0$.
Hence,  by \eqref{lambda_zero},  $\varphi_{2}$ never vanishes on $I$. 
By \eqref{segno}, $u_2$ is constantly equal to $1$ or $-1$ on $I$.
From the first equation in \eqref{eq:pmpphi} we have $\varphi_3=0$ on $I$, and hence on the whole interval of definition of the trajectory.
In particular, by  \eqref{eq:pmpphi} we have that the whole trajectory is $\varphi_1$-singular.

Conversely,   every trajectory   corresponding to 
$u_2=\pm 1$ constant and $u_1$ measurable with $|u_{1}|\leq 1$ has a  $\varphi_1$-singular extremal lift 
with $\varphi_2=1$ and $\varphi_1=\varphi_3=0$. 

Moreover, each of such curves $\bar \gamma = (\bar x, \bar y , \bar z) : [0,T] \to \mb{H} $ is 
time-minimizing since
$T=|\bar y(0) - \bar y(T)|$
and $|\dot y| \leq 1$ for every trajectory of $\dot \gamma= u_1 X_1(\gamma) + u_2 X_2(\gamma)$, with $|u_1|,|u_2|\leq 1$.
\end{proof}

\subsection{Regular arcs}
\bl \label{lemma-bang-H}
 On $\mb{H}$
the   trajectories that have   a regular arc are regular bang-bang. Moreover, all arcs have the same length $s$ except possibly the last and the first arc, whose lengths are less than or equal to $ s$. 
At the junction between regular arcs the components  $u_1$ and $u_2$ of the control   switch sign alternately. 
\el


%
%
%

\begin{figure}[h]
\begin{center}
	\scalebox{.75}{\input{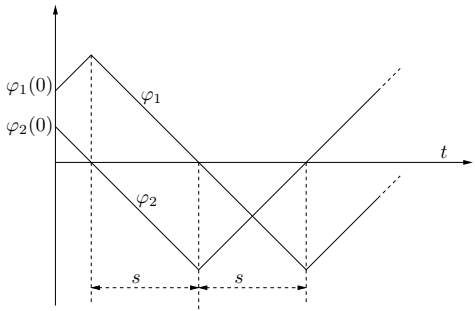}}
	\caption{The switching functions for the Heisenberg structure, when $\ph_{3}\neq0$.}
\label{figura1}
\end{center}
\end{figure}

\begin{proof}
Let $I$ be an interval corresponding to a regular arc of the trajectory.
Without loss of generality, 
$\varphi_{1},\varphi_{2}>0$  on $I$. Hence, by \eqref{segno} we have
$u\equiv (1,1)$ on $I$.
Fix $ t_0\in I$.
Two cases are possible:
\begin{description}
\iii[(a) $\varphi_{3}(t_0)=0$]
By  \eqref{eq:pmpphi} we have that $\varphi_1$ and  $\varphi_2$ are constant along the entire trajectory, which is then a single regular arc.

\iii[(b) $\varphi_{3}(t_0)\neq 0$]
Denote by $a$ the constant value of $\varphi_{3}$.
Using \eqref{eq:pmpphi} we find
$$\varphi_{1}(t)=\varphi_{1}(t_0)-a(t-t_0),\qquad \varphi_{2}(t)=\varphi_{2}(t_0)+a(t-t_0), \qquad \forall t\in I.$$
\end{description}
Without loss of generality $a>0$. Set $t_1=t_0+\varphi_1(t_0)/a$.
If the trajectory is defined up to time $t_1$, then  $\varphi_1$ and  $\varphi_2$ are positive in the interval $(t_0,t_1) $.
Also $\varphi_1(t_1)=0$.

Since $u_2=1$ in a neighborhood of $t_1$, we deduce that $\varphi_1$ is affine in a neighborhood of $t_1$, with slope $-a$.  Hence $\varphi_1<0<\varphi_2 $  in a right-neighborhood of $t_1$. 
Then $t_1$ is the starting time of another regular arc with control  $u=(-1, 1)$.

Repeating this argument, backwards in time as well, we conclude that the extremal trajectory is the concatenation of regular arcs
of length
  $\varphi_2(t_1)/a=(\varphi_1(t_0)+\varphi_2(t_0))/\varphi_3(t_0)$, except possibly for the first and last arc, see Figure \ref{figura1}.  
  The switching occur alternately for $u_1$ and $u_2$. 
\end{proof}

The picture of the switching function, in the nontrivial case (b), is given in Figure \ref{figura1}.

\subsection{Bound on number of optimal regular arcs}

%
%

\bp \label{lemma-bang-optimal-H}
A regular bang-bang trajectory with more than $5$  arcs is not optimal.
\ep

\begin{proof}


Let us consider a trajectory  with $6$ bang arcs. 
By Lemma \ref{lemma-bang-H}, without loss of generality
 we can assume that  the  
successive values of the control  are 
$$(-1,-1), (-1,1), (1,1), (1,-1), (-1,-1), (-1,1).$$
Denote the length of the internal bang arcs by $s$ (recall that all arcs have the same, except possibly the first and last).

We are going to apply Theorem \ref{thm2nd} by taking $j=3$. 
Let
$\tau_3$ be the third switching time.
Since at $\tau_3$ the function $\varphi_2$  switches  sign, we have that
 $\varphi_2(\tau_3)   = \la \lambda(\tau_3),  X_2(\gamma(\tau_3))\ra=0$. 
 
 Up to multiplication of  $\lambda(\cdot)$ by a positive scalar, we can 
normalize $\varphi_3$, which is constant, to $-1$.
Hence,
$\varphi_1(\tau_3)=s$, which implies that
 $\lambda(\cdot)$ is uniquely determined by the sequence of switching times.
Set
 $$X_+=X_1+ X_2,\quad X_-=X_1-X_2.$$
Following the notations of Theorem \ref{thm2nd}, we have
\begin{align*}
Z_0&= e^{-s\,\ad (X_+)} e^{s\,\ad (X_-) } (-X_+)=-X_+-2s X_{3},\\
Z_1&= e^{-s\,\ad (X_+)}(-X_-)=-X_--2s X_{3},\\
Z_2&=X_+,\\
Z_3&=X_-,\\
Z_4&=e^{s\,\ad (X_-)}(-X_+)=-X_+-2s X_{3},\\
Z_5&=e^{s\,\ad (X_-)}e^{s\,\ad (-X_+)}(-X_-)=-X_--2 s X_{3}.
\end{align*}
A simple calculation shows that
\begin{align*}		
\sigma_{01}&=\sigma_{05}=\sigma_{12}=\sigma_{23}=\sigma_{34}=\sigma_{45}=2,\\
\sigma_{02}& =\sigma_{04}=\sigma_{13}=\sigma_{15} =\sigma_{24} =\sigma_{35}  =0,\\
\sigma_{03}&=\sigma_{14}=\sigma_{25}=-2.
\end{align*}
Decomposing the relation $\sum_{i=0}^5 \alpha_i Z_i(\gamma(\tau_3))=0$ on the basis $\{X_+(\gamma(\tau_3)),X_-(\gamma(\tau_3)),X_{3}(\gamma(\tau_3))\}$, one gets
\begin{align*}
-\alpha_0+\alpha_2-\alpha_4=0,\quad
-\alpha_1+\alpha_3 -\alpha_5=0,\quad 
2s (-\alpha_0-\alpha_1-\alpha_4-\alpha_5)=0.
\end{align*}
Solving in $\alpha_0,\alpha_1,\alpha_2$, gives
\begin{align*}
\alpha_3= -\alpha_2,\quad
\alpha_4= -\alpha_0 + \alpha_2,\quad
\alpha_5=-\alpha_1-  \alpha_2.
\end{align*}
Notice that the relation $\sum_{i=0}^5 \alpha_i=0$ is automatically satisfied. 
Then we can parameterize the  space $W$ appearing in the statement of Theorem \ref{thm2nd} by 
$\alpha_0,\alpha_1,\alpha_2$, i.e.,  
$$W=\{(\al_0,\al_1,\al_2,-\al_2,-\al_0+\al_2,-\al_1-\al_2)\mid \al_0,\al_1,\al_2\in\R\},$$ and write the quadratic form $Q$ as 
$$Q(\alpha_0,\alpha_1,\alpha_2)=  4 \alpha_0\alpha_1 + 4 \alpha_0 \alpha_2  -4 \alpha_2^2.$$
In particular,
$Q(1,1,0)=4>0$, which implies that the trajectory is not optimal.
\end{proof}

\subsection{Optimal trajectories and shape of the unit ball}
Here we summarize the results obtained in the previous sections and we plot the unit ball in the Heisenberg group.

Recall that once we characterize the controls $u_{1}(t)$ and $u_{2}(t)$ associated with an extremal trajectory,  the trajectory itself can be recovered by solving the differential equation
$$\dot \g(t)=u_{1}(t)X_{1}(\g(t))+u_{2}(t)X_{2}(\g(t)).$$
By the coordinate expression \eqref{eq:campiH} of the vector fields $X_{1},X_{2}$, this is equivalent to solving the system
\bqn
\begin{cases}
\dot x=u_{1}\\
\dot y=u_{2}\\
\dot z=\frac12 (u_{2}x-u_{1}y)
\end{cases}
\eqn 
In particular, the trajectory is determined by its  projection $\tilde\g$ onto the $xy$-plane, since the $z$ coordinate of the trajectory can be found by integration. As it is well-known, it computes the signed area defined by the closed curve given by following $\tilde\g$ and then coming back to the origin along a line segment.

\begin{figure}
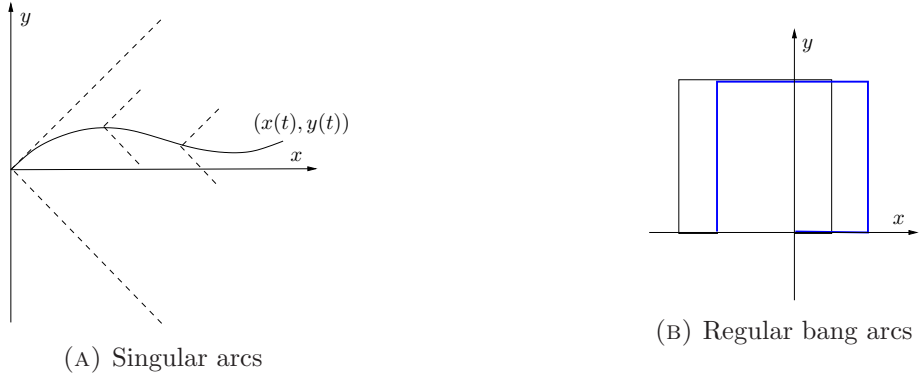

\begin{subfigure}{.5\textwidth}
  \centering
 \scalebox{.8}{\input{figura-2a.pstex_t}}
	\caption{Singular arcs}
\label{figura2a}	
\end{subfigure}%
\begin{subfigure}{.5\textwidth}
  \centering
  \scalebox{.8}{\input{figura-2b.pstex_t}}
	\caption{Regular bang arcs}
  \label{figura2b}
  
\end{subfigure}
\caption{Singular and regular arcs in the Heisenberg group}
\end{figure}

As discussed in Lemma \ref{lemma-singular-H}, the singular trajectories correspond 
to the case when the control $u_{1}(t)$ is constantly equal to $\pm 1$ and $u_{2}(t)$ is free (or the symmetric situation). In Figure \ref{figura2a} we can see an example of such a curve when $u_{1}(t)=1$. Recall that these curves are optimal for all times and that, given one such trajectory,  there exists a time-optimal bang-bang trajectory with at most 3 bang arcs connecting the same endpoints.


Regular bang-bang trajectories correspond to switching functions as in Figure~\ref{figura1}, where the controls switch sign alternately. These trajectories draw squares in the $xy$-plane as in Figure~\ref{figura2b}.  

If such a trajectory has more than 5 bang arcs, then Proposition \ref{lemma-bang-optimal-H} guarantees that the trajectory is not optimal.

Notice that there exist time-minimizing curves of this kind with
$5$ regular bang arcs, as illustrated in Figure~\ref{figura2b}. However, not all bang-bang trajectories with 
$5$ bang arcs are time-minimizing. Indeed, if the underlying square is swept more than once, then the trajectory is no more a minimizer. Finally, let us also remark that for every minimizer with 5 regular bang arcs there exists a minimizer with 4 regular bang arcs joining the same endpoints (see again Figure~\ref{figura2b}). 

\begin{figure}
\begin{subfigure}{.5\textwidth}
  \centering
\scalebox{0.18}{\includegraphics{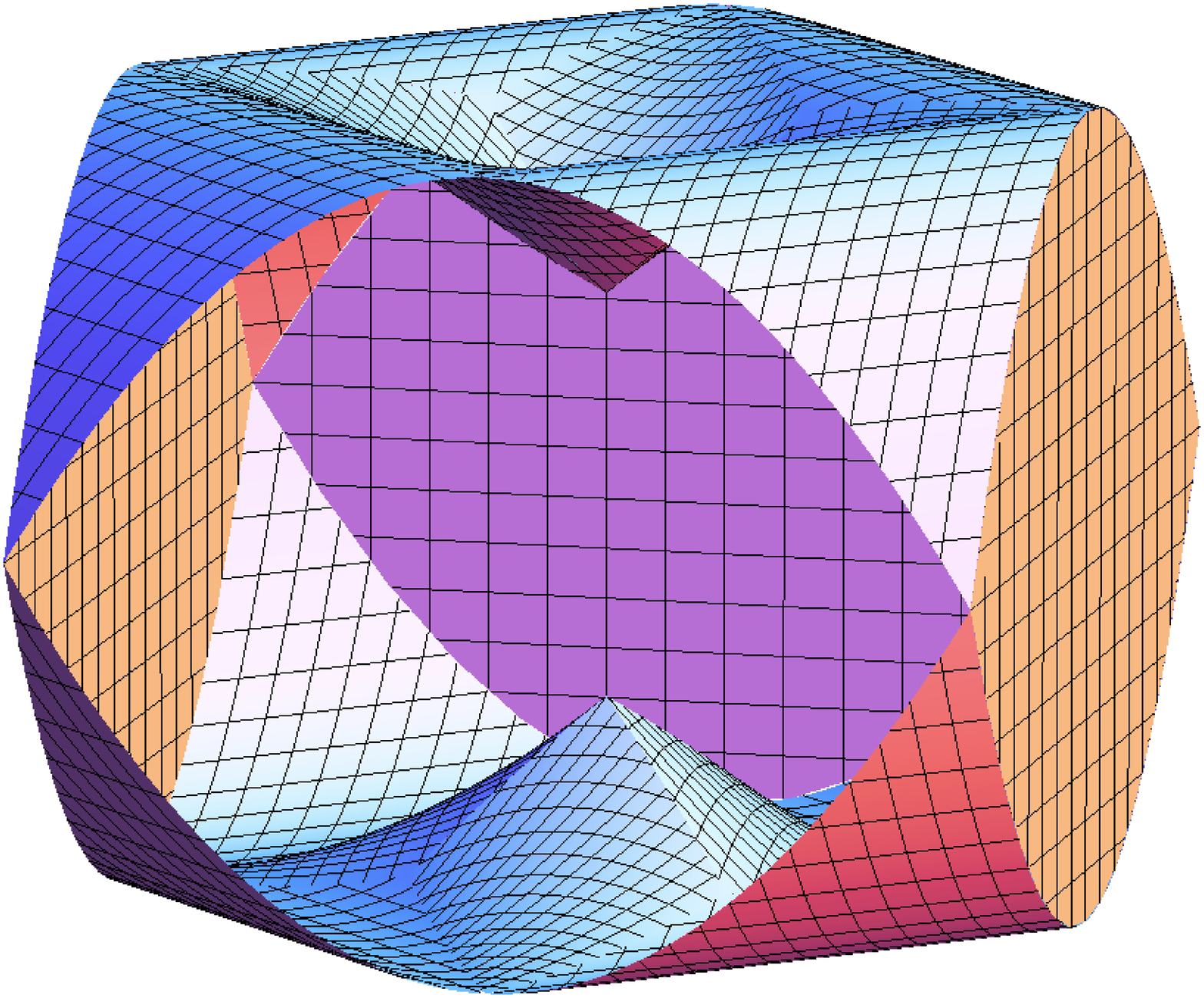}}
	\caption{Sphere}
\label{figura3H}	
\end{subfigure}%
\begin{subfigure}{.5\textwidth}
  \centering
  \scalebox{0.18}{\includegraphics{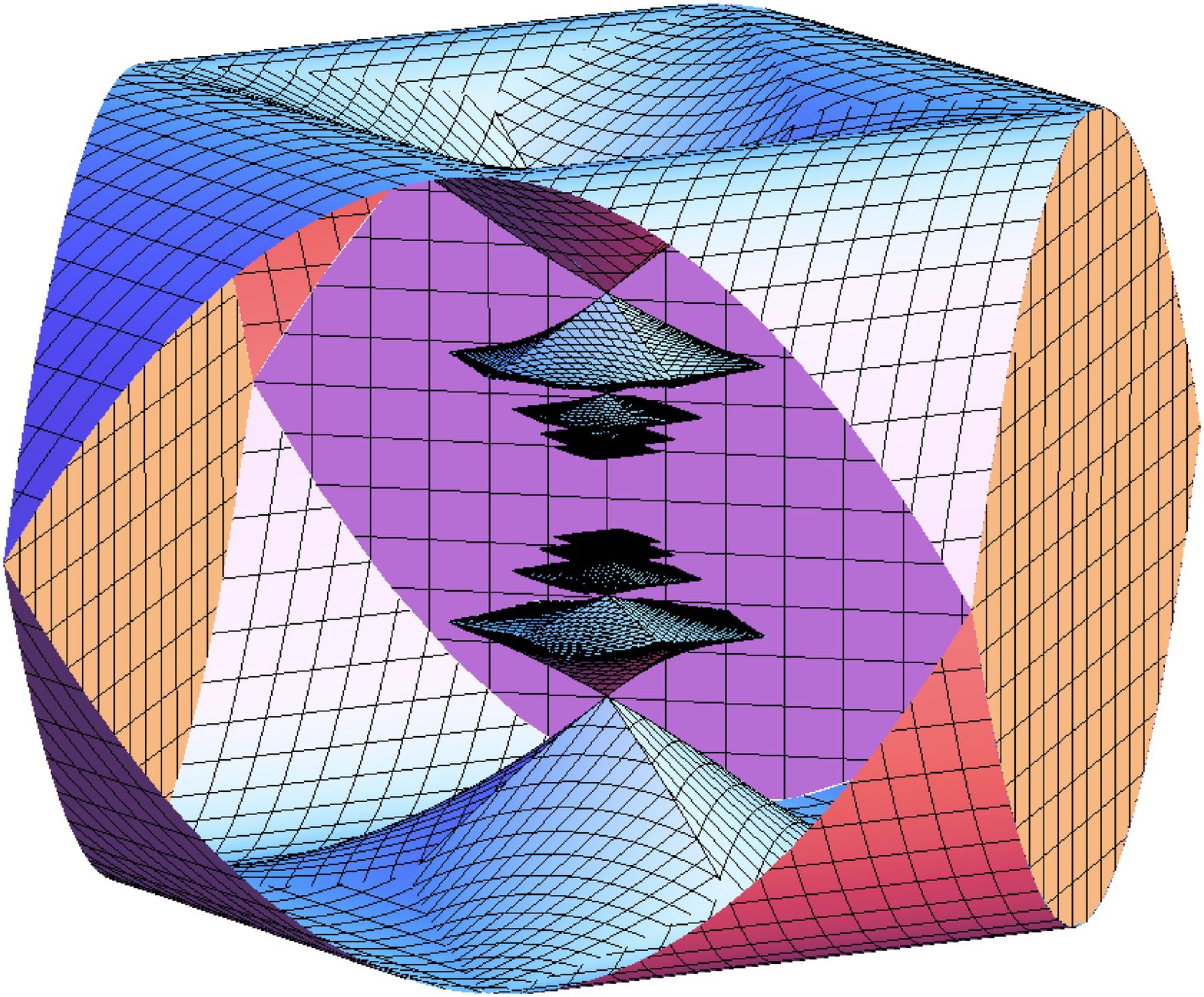}}
	\caption{Front}
  \label{figura3bH}
  
\end{subfigure}
\caption{Sphere and front of the unit sphere in the Heisenberg group}
\end{figure}

We stress that, by the classification of the previous sections, the only extremal trajectories connecting two distinct points on the same vertical line in the Heisenberg group are regular bang-bang. 
Once the shape of optimal trajectories is known, a picture of the Heisenberg sphere can be easily drawn. See Figure~\ref{figura3H} and Figure~\ref{figura3bH}.

\section{Grushin structures}
In this section we provide a description of the time-minimizing trajectories in two different sub-$\ell^\infty$ structures in the Grushin plane.

The classical sub-Riemannian structure on the Grushin plane is the metric structure on $\R^{2}$ determined by the choice of the orthonormal vector fields 
$$X_{1}=\partial_{x} ,\qquad X_{2}=x \partial_{y} .$$
In other words the sub-Riemannian distance is characterized as follows
$$d(p_{1},p_{2})=\inf \left\{\int_{0}^{T} \sqrt{u_{1}^{2} +u_{2}^2}dt\, \bigg|\, \dot \gamma=u_{1}X_{1}(\gamma)+u_{2}(t)X_{2}(\gamma),\, \gamma(0)=p_{0},\, \gamma(1)=p_{1}  \right\}.$$
The geodesic problem for this distance is equivalent to the time-optimal control problem defined by $X_{1}$, $X_{2}$ and $u(t)\in B$, where $B=\{u_{1}^{2} +u_{2}^2\leq 1\}$ is the standard Euclidean ball.\\

Due to the lack of symmetry of the Grushin structure, it is meaningful to consider two different   sub-$\ell^\infty$ structures on the 
Grushin plane.

\subsection{The first structure}
Consider the sub-$\ell^\infty$ structures on
 $  \R^2$
determined by the vector fields 
\bqn\label{eq:gru0}
Y_{1}=\partial_{x}+x \partial_{y} ,\qquad Y_{2}=\partial_{x}- x \partial_{y} .
\eqn
Notice that $Y_{1}=X_{1}+X_{2}$ and $Y_{2}=X_{1}-X_{2}$, so that we are considering the sub-$\ell^{1}$ sub-Finsler structure associated with $X_1$, $X_2$, up to a dilation factor.
Similarly as in the Heisenberg group, let us introduce the vector field $Y_3=\partial_{y}$.

The Lie algebra generated by $Y_{1},Y_{2},Y_{3}$ actually satisfies the same commutator relations as in the Heisenberg group, namely 
\bqn\label{derivata-switch0}
[Y_1,Y_2]=Y_3,\qquad [Y_1,Y_3]=[Y_2,Y_3]=0.
\eqn
The identities \eqref{derivata-switch0} gives the same equations \eqref{eq:pmpphi} obtained in the Heisenberg case for the switching functions along an extremal trajectory
\bqn \label{eq:pmpphigru0}
\dot{\varphi}_{1}=-u_{2}\varphi_{3},\qquad \dot{\varphi}_{2}=u_{1}\varphi_{3}, \qquad \dot{\varphi}_{3}=0.
\eqn
In particular, $\varphi_3$ is constant.
From the relation $Y_2-Y_1=x Y_3$ we  have the additional relation
\bqn\label{eq:gru01}
\varphi_2 - \varphi_1 = x \varphi_3.
\eqn
In the case $\varphi_3=0$,
we get $\varphi_2  \equiv \varphi_1$ equals  a constant, which is different from zero since the covector cannot be identically zero. Hence $u_1,u_2$ are both $1$ or both $-1$.
In other words, every trajectory corresponding to $\varphi_3=0$ is a horizontal line. Such curves are indeed time-minimizers.



Let us then consider the case $\varphi_3\neq0$. Under this assumption, if both $\varphi_{1}$ and $\varphi_{2}$ are vanishing then the trajectory is abnormal. In particular relation \eqref{eq:gru01} implies that $x(t)=0$ along the trajectory, hence we deduce that the trajectory is reduced to a point which is contained in the $y$-axis.

\bl \label{lemma-abnormal-gru0}
 The only abnormal arcs on the sub-$\ell^{\infty}$ structure on $\R^{2}$ defined by the vector fields \eqref{eq:gru0} are the constant curves contained in the set $\{x=0\}$.
 Consequently, no minimizer joining two distinct points is abnormal.
\el

An analogous reasoning shows that there are no $\varphi_{1}$-singular (resp.\ $\varphi_{2}$-singular) trajectories that are not abnormal. Indeed assume the trajectory is $\varphi_{1}$-singular. Then $\varphi_{1}(t) =0$ for all $t$, that implies, by \eqref{eq:pmpphigru0}, that $u_{2}=0$ for all $t$ (recall that we are in the case $\varphi_{3}\neq0$). Thus $\varphi_{2}$ is 
necessarily identically zero and the trajectory is actually abnormal. The situation is analogue for $\varphi_{2}$-singular trajectories.

Following the lines of  Lemma \ref{lemma-bang-H} one can then show that, if the trajectory is not abnormal, then it is bang-bang, and all internal arcs of a bang-bang trajectory have the same length $s$. 

\bl \label{lemma-bang-gru0}
 On  the sub-$\ell^{\infty}$ structure on $\R^{2}$ defined by the vector fields \eqref{eq:gru0}
the   trajectories that have   a regular arc are regular bang-bang. Moreover, all arcs have the same length $s$ except possibly the last and the first arc, whose lengths are less than or equal to $ s$. 
At the junction between regular arcs the components  $u_1$ and $u_2$ of the control   switch sign alternately. 
\el

\subsubsection{Bound on number of optimal regular arcs}
Notice that  if
a bang-bang trajectory
has an internal bang arc whose length is $t$, then
$u_1$  and $u_2$ switch on the lines  $x=\pm t/2$ (see Figura~\ref{figura4}).

Regarding optimality, we claim that a bang-bang trajectory with   $4$ bang arcs is not optimal.
Indeed,  every extremal trajectory starting from $y$-axis is not optimal after it intersects again  the vertical axis, as it follows by replacing the trajectory by its reflection along the $y$-axis.

 \subsubsection{Optimal trajectories and shape of the unit ball}

\begin{figure}[h]
\begin{center}
	\scalebox{.9}{\input{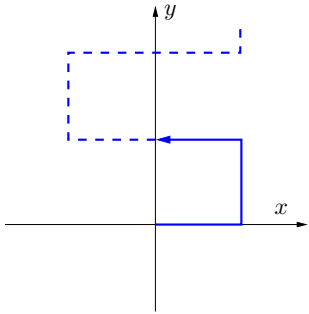}}
	\caption{Regular bang-bang trajectories for the Grushin structure \eqref{eq:gru0}}
\label{figura4}
\end{center}
\end{figure}

\begin{figure}[h]
\begin{center}
\includegraphics[scale=0.6]{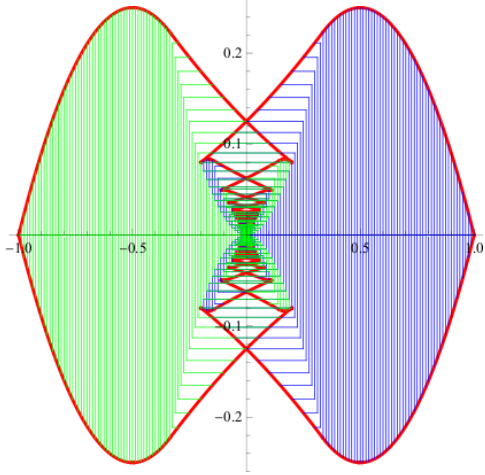}
	\caption{The unit sphere and its front for the Grushin structure \eqref{eq:gru0}}
\label{figura5}
\end{center}
\end{figure}

The picture of the regular bang bang trajectories for this structure on the Grushin plane is given in Figure \ref{figura4}. The corresponding picture of the unit ball is obtained in Figure \ref{figura5}.

\subsection{The second structure}
We consider now the sub-$\ell^\infty$ structure on
 $  \R^2$
determined by the vector fields 
\bqn\label{eq:gru00}
X_{1}=\partial_{x} ,\qquad X_{2}=x \partial_{y} .
\eqn
and we introduce the vector field $X_3=\partial_{y}$.

The Lie algebra generated by $X_{1},X_{2},X_{3}$ again satisfies the same commutator relations as in the Heisenberg group, namely 
$$[X_1,X_2]=X_3,\qquad [X_1,X_3]=[X_2,X_3]=0.$$
Thus the identity \eqref{derivata-switch} gives the analog equations \eqref{eq:pmpphi} for the switching functions along an extremal trajectory
\bqn \label{eq:pmpphigru00}
\dot{\varphi}_{1}=-u_{2}\varphi_{3},\qquad \dot{\varphi}_{2}=u_{1}\varphi_{3}, \qquad \dot{\varphi}_{3}=0,
\eqn
In particular, $\varphi_3$ is constant. 
From $X_2=x X_3$ we  have the additional relation
$$\varphi_2 = x \varphi_3.$$

In the case $\varphi_3=0$,
we get $\varphi_2  \equiv 0$ and $\varphi_1$ equals  a nonzero constant (otherwise the covector is identically zero).
Reasoning as in Lemma \ref{lemma-singular-H}, 
we have immediately the following result

\bl \label{lemma-singular-gru1}
 On the sub-$\ell^{\infty}$ structure on $\R^{2}$ defined by the vector fields \eqref{eq:gru00}
the nonconstant trajectories that have   singular arcs are 
exactly those
for which
$u_1$ is constantly equal to $1$ or $-1$.
All of them consist of a single singular arc and are time-minimizers.
\el

Let us then assume in what follows $\varphi_3\neq0$.

\bl \label{lemma-abnormal-gru00}
 The only abnormal arcs on the sub-$\ell^{\infty}$ structure on $\R^{2}$ defined by the vector fields \eqref{eq:gru00} are the constant curves contained in the $y$-axis. 
 Consequently, no minimizer joining two distinct points is abnormal.
\el
If $\varphi_3\neq0$ and the trajectory is not abnormal, then 
as in Lemma \ref{lemma-bang-H}
it is regular bang-bang, all arcs  have the same length $s$ except possibly the last and the first one, whose lengths are less than or equal to $ s$. 
At the junction between bang arcs the components  $u_1$ and $u_2$ of the control   switch sign alternately. 

Moreover, on a regular bang-bang trajectory,
$u_2$    switches on the line  $x=0$, since, if $\varphi_2(t)=0$ at a point $t$, then $x(t) \varphi_3 = 0$.
Therefore if
a bang-bang trajectory
has an internal bang arc whose length is $s$, then
$u_1$   switches on the lines  $x=\pm s$. Moreover, 
at $u_1$-switching times the function $u_1$ goes from $1$ to $-1$ if the switch occurs in the half-plane $x>0$
while it goes from $-1$ to $1$   in the half-plane $x<0$, since 
$$ {\rm sign}(\dot\varphi_{1})= -  {\rm sign}( u_{2}\varphi_{3} ) =- {\rm sign}( \varphi_{2}\varphi_{3} ) 
=- {\rm sign}( x \varphi_{3}^2 ) =- {\rm sign} (x) .$$

\subsubsection{Bound on number of optimal regular arcs}

Regarding optimality, we  prove the following lemma.
\bl \label{lemma-bang-optimal-grushin}
A regular bang-bang trajectory with more than $3$  arcs is not optimal.
If, moreover, the trajectory starts on the $y$-axis and it is optimal, then it has at most $2$ arcs. 
\el
\begin{proof}
First notice that, contrarily to what happens in the Heisenberg case, the role of the two vector fields $X_1$, $X_2$ is not symmetric. The replacement of 
$(u_1,u_2)$ by $(-u_1,-u_2)$ coupled with the reversion in the order of bangs, on the contrary, still yields a symmetric, equivalent, situation. This is a general fact, since it simply corresponds to reverse the parameterization of the curve.
Looking at regular bang-bang trajectories (see Figure~\ref{figura7}) one 
immediately recognizes that 
the proof of the lemma can be given by looking at
two types of bang-bang trajectories, whose successive values of the control are
$$ (1,-1), (1,1), (-1,1), (-1,-1)\quad \mbox{ and }\quad (1,1), (-1,1), (-1,-1), (1,-1),$$ 
respectively.
In the first case, one notices that reflecting the second and third bang arcs with respect to the $y$-axis
yields another horizontal curve with the same length, which is not extremal. 
Hence the curve is not optimal. This argument also shows that 
regular bang-bang trajectories starting from the $y$-axis  and with more than 2 bang arcs are not optimal. 

In the second case, let us apply Theorem \ref{thm2nd} at the second switching time. One gets
$$
Z_0  =X_1+X_2 + 2s X_{3},\qquad
Z_1=-X_1+X_2,\qquad
Z_2=-X_1-X_2,\qquad
Z_3= X_1-X_2 +2s X_{3}.
$$
Parameterizing the  space $W$
by the coordinates $\alpha_0,\alpha_1$ we get that $W=\{(\al_0,\al_1,-\al_1,-\al_0)\mid \al_0,\al_1\in \R\}$. 
Normalizing $\varphi_3=1$ (uniqueness of the covector up to a positive factor is proved as in the case of the Heisenberg group),  we write the quadratic form $Q$ as 
$$Q(\alpha_0,\alpha_1)= 2 \alpha_0 ^2 +4  \alpha_0 \alpha_1 -2  \alpha_1^2.$$
Since $Q(1,0)$ is positive,  the considered trajectory 
is not optimal.
This concludes the proof of Lemma~\ref{lemma-bang-optimal-grushin}.
\end{proof}

 \subsubsection{Optimal trajectories and shape of the unit ball}
 
\begin{figure}[h]
\begin{center}
	\scalebox{.8}{\input{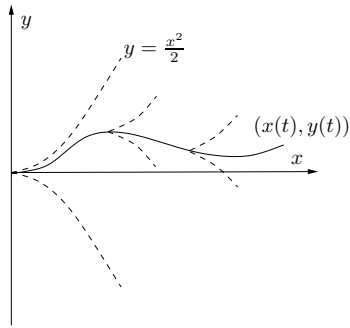}}
	\caption{Singular trajectories for the Grushin structure \eqref{eq:gru00}}
\label{figura6}
\end{center}
\end{figure}

\begin{figure}[h]
\begin{center}
	\scalebox{.9}{\input{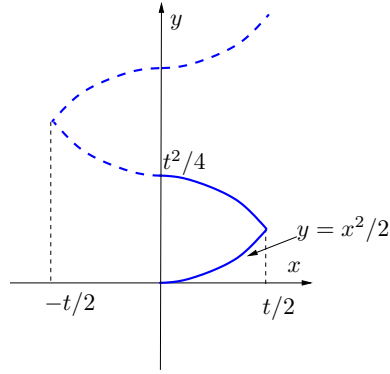}}
	\caption{Bang-bang trajectories for the Grushin structure \eqref{eq:gru00}}
\label{figura7}
\end{center}
\end{figure}

\begin{figure}
\begin{subfigure}{.5\textwidth}
  \centering
\scalebox{0.14}{\includegraphics{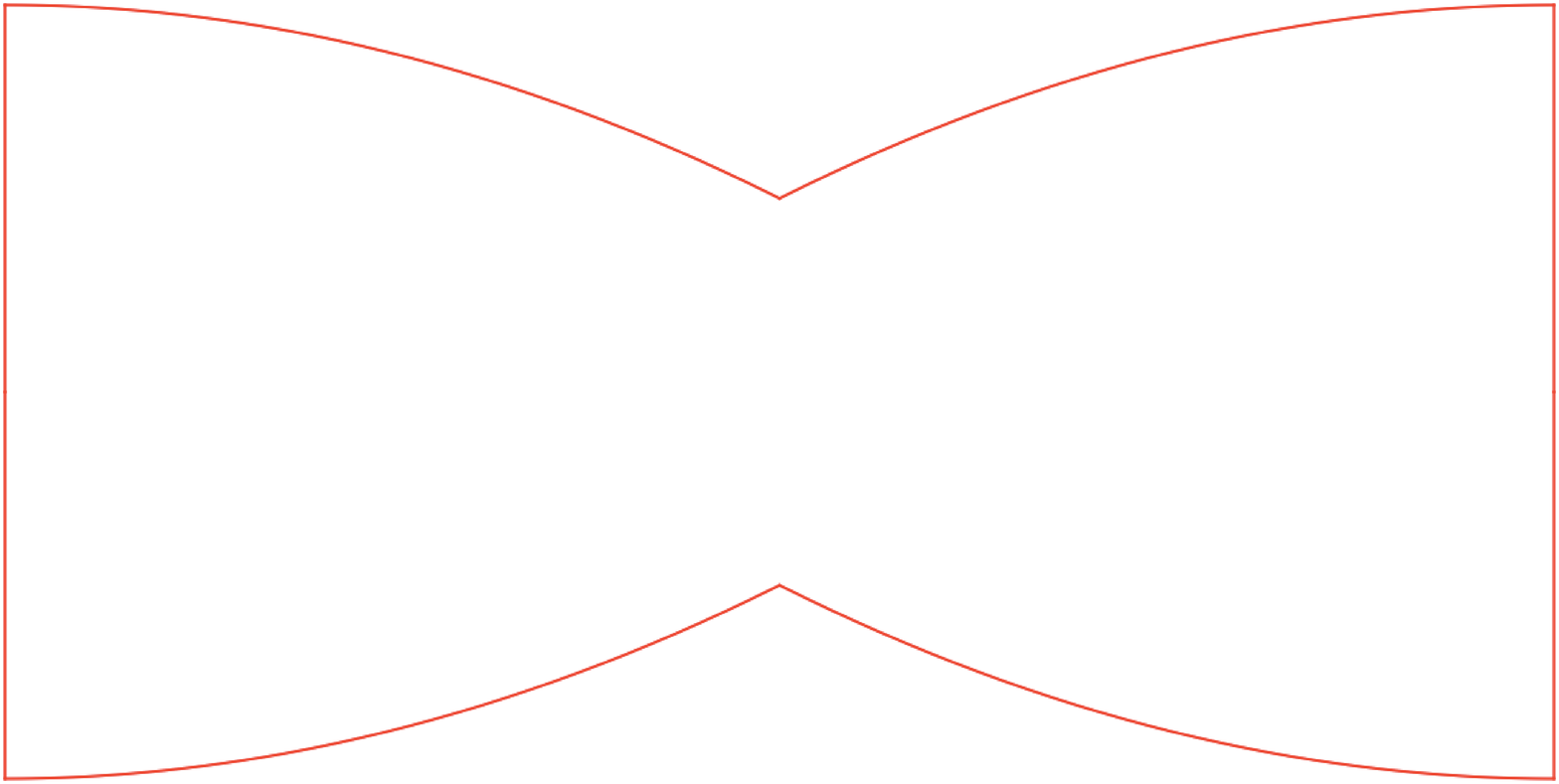}}
	\caption{Sphere}
\label{figura8a}	
\end{subfigure}%
\begin{subfigure}{.5\textwidth}
  \centering
  \scalebox{0.14}{\includegraphics{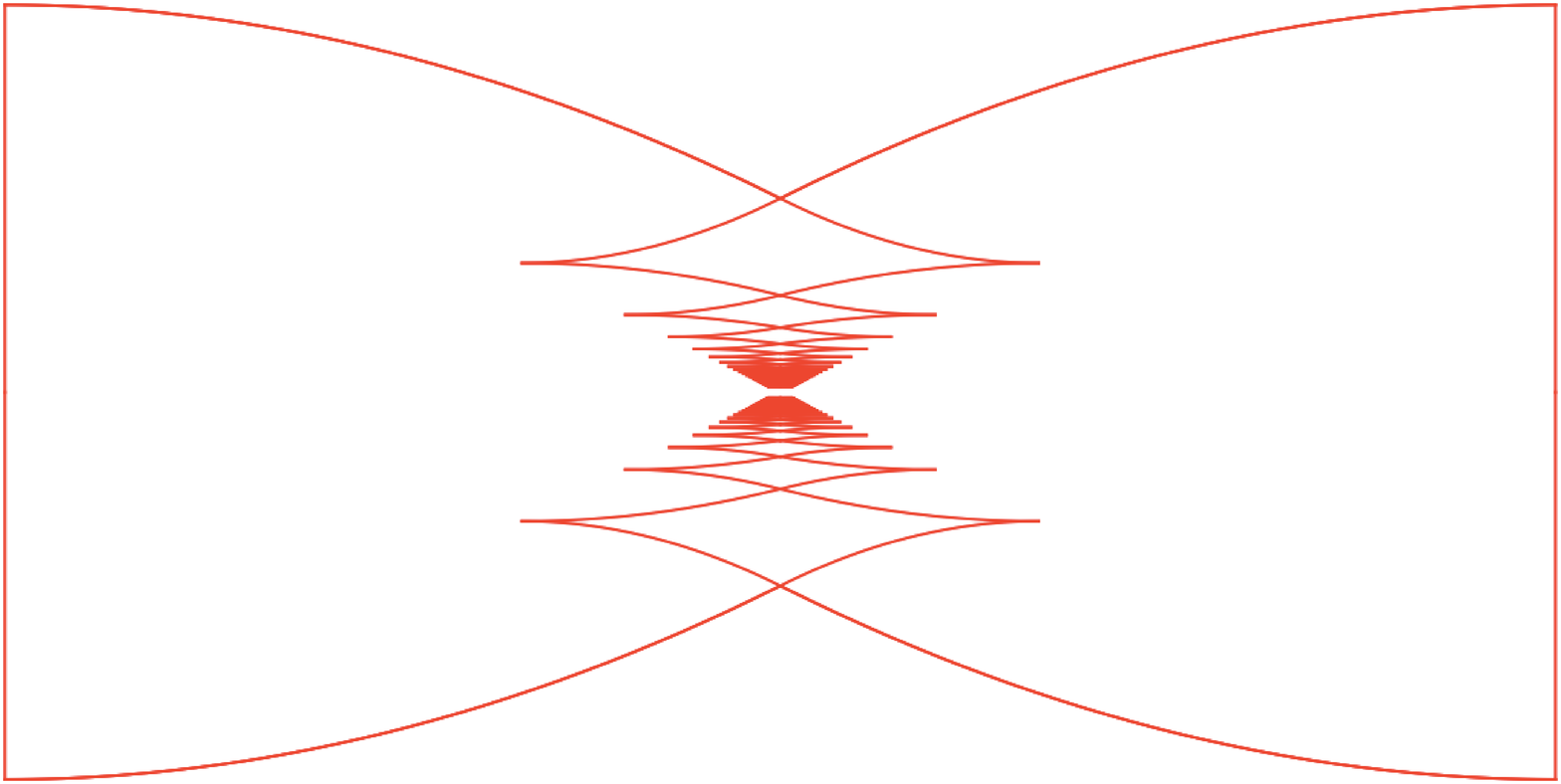}}
	\caption{Front}
  \label{figura8b}
  
\end{subfigure}
\caption{Sphere and front of the unit sphere for the Grushin structure \eqref{eq:gru00} \label{figura8}}
\end{figure}

 In this structure for the Grushin plane we have singular trajectories that are similar to the one obtained in the Heisenberg group, see Figure \ref{figura6}. Let us stress that in this case the tangent vector of the curve is forced to be inside a cone whose width increases with the $x$ coordinate.
Regular bang bang trajectories from the origin are illustrated in Figure~\ref{figura7}. These trajectories lose optimality as soon as they reach the vertical axes.
The picture of the unit ball in the Grushin plane with this structure is in Figure \ref{figura8}.

\section{Martinet structures}\label{s:martinet}
In this section we provide a description of the time-minimizing trajectories for two different sub-$\ell^\infty$ structures associated with the Martinet distribution. This is the easiest example where nontrivial abnormal minimizers appear.

The classical sub-Riemannian structure on the Martinet space is the metric structure on $\R^{3}$ determined by the choice of the orthonormal vector fields 
$$X_{1}=\partial_{x}+y^{2} \partial_{z},\qquad X_{2}=\partial_{y}.$$
The sub-Riemannian distance is then 
$$d(p_{1},p_{2})=\inf \left\{\int_{0}^{T} \sqrt{u_{1}^{2} +u_{2}^2}dt\, \bigg|\, \dot \gamma=u_{1}X_{1}(\gamma)+u_{2}(t)X_{2}(\gamma),\, \gamma(0)=p_{0},\, \gamma(1)=p_{1}  \right\}.$$
The geodesic problem for this distance is equivalent to the time-optimal control problem defined by $X_{1},X_{2}$ and $u(t)\in B=\{u_{1}^{2} +u_{2}^2\leq 1\}$.

As in the case of the Grushin plane, due to the lack of symmetry, we are lead to consider two different   sub-$\ell^\infty$ structures. 

\subsection{The first structure} \label{s:mart0}
Consider the sub-$\ell^\infty$ structure on
 $  \R^3$
determined by the vector fields 
\bqn\label{eq:mart0}
Y_{1}=\partial_{x}+\partial_{y}+y^{2}\partial_{z},\qquad Y_{2}=\partial_{x}-\partial_{y}+y^{2}\partial_{z}.
\eqn
Notice that $Y_{1}=X_{1}+X_{2}$ and $Y_{2}=X_{1}-X_{2}$, so that we are considering the sub-$\ell^{1}$ Finsler structure defined by $X_{1}$ and $X_{2}$, up to a dilation factor.
In analogy to the other cases, let us introduce the following vector fields defined by the commutators of the elements of the basis of the distribution
\begin{equation}\label{BrRe}
Y_{3}:=[Y_{1},Y_{2}]=4y\partial_{z},\quad Y_{4}:=[Y_{1},[Y_{1},Y_{2}]]=4\partial_{z},\quad Y_{5}:=[Y_{2},[Y_{1},Y_{2}]]=-4\partial_{z}.
\end{equation}
The switching functions associated with these vector fields and with an extremal pair $(\lambda(\cdot),\gamma(\cdot))$ are
$$\varphi_{i}(t)=\la \lambda(t),Y_{i}(\gamma(t))\ra,\qquad i=1,\ldots,5.$$
They 
satisfy the following system of differential equations
\begin{gather} \label{eq:pmpphimart0}
\dot{\ph}_{1}=-u_{2}\ph_{3},\qquad \dot{\ph}_{2}=u_{1}\ph_{3}, \qquad \dot{\ph}_{3}=u_{1}\ph_{4}+u_{2}\ph_{5},\\
\dot{\ph}_{4}=0,\qquad \dot{\ph}_{5}=0. \nn
\end{gather}
\brem \label{r:ph4} 
It follows from the bracket relations \eqref{BrRe}  that $\ph_{4}$ and $\ph_{5}=-\ph_{4}$ are constants 
and we have $\ph_{3}=y\ph_{4}=-y\ph_{5}$. In particular, if $\ph_4=0$, then $\ph_3$ is also constantly equal to zero, and $\ph_1,\ph_2$ are constant.
%
\erem
\bl \label{l:reasoning}
The nontrivial abnormal arcs on the sub-$\ell^{\infty}$ structure on $\R^{3}$ defined by the vector fields \eqref{eq:mart0} are the horizontal lines contained in the plane $\{y=0\}$.
\el
\begin{proof} 
Assume that the trajectory is not reduced to a point and it is abnormal on some interval $I$. In particular we have $\ph_{1}(t)=\ph_{2}(t)=0$ for all $t\in I$, while its control $(u_{1}(t),u_{2}(t))$ is not identically zero on $I$.
From identities \eqref{eq:pmpphimart0} one immediately gets that $-u_{2}(t)\ph_{3}(t)=u_{1}(t)\ph_{3}(t)=0$. Hence, if  we have that $\ph_{3}(t)=y(t) \ph_{4} =0$ for every $t$ (recall that $\ph_{4}$ is constant), 
then $y(t)=0$ for all $t\in I$, otherwise $\ph_{4}=0$ and the covector is identically zero. 
In particular $u_1=u_2$ on $I$ and  the trajectory is contained in a line $\{y=0,z=z_{0}\}$.
\end{proof}
\brem Every costant trajectory (i.e., such that $u_{1}(t)=u_{2}(t)=0$ on $[0,T]$) is also abnormal.
\erem

The classification of extremal trajectories on Martinet is then reduced to regular and those that are  singular with respect to exactly one control.

\subsubsection{Singular arcs}
Let us now consider a singular arc. We show that in this case we can recover its (singular) control by differentiation of the adjoint equations.

Indeed assume that the trajectory is $\ph_{1}$-singular, i.e., $\ph_{1}\equiv0$ on $I$, and we want to recover its associated control $u_{1}$. Notice that $|u_{2}|=1$ is constant and $\dot \ph_{1}=-u_{2}\ph_{3}$. 
By singularity assumptions $\dot{\ph}_{1}(t)\equiv 0$, that implies $\ph_{3}(t)\equiv 0$ for all $t\in I$. 
We deduce that either $\ph_{4}=0$ or $u_{1}=u_{2}$ on $I$. We have two possibilities::
\bi
\iii[(i)] if $\ph_{4}=0$ then $u_{1}$ is free, 
\iii[(ii)] if $\ph_{4}\neq0$ then $u_{1}=u_{2}$ and the singular arc is also a bang arc, 
with no constraint on its length. Moreover $y=0$ on such an arc. 
\ei
The situation with $\ph_{2}$-singular arcs is perfectly symmetric.
\subsubsection{Regular arcs}
Assume that both $\ph_{1}(0),\ph_{2}(0)\neq 0$ (the same holds for small times by continuity).
Because of Remark~\ref{r:ph4}, we can assume that $\ph_4\ne 0$ (otherwise the trajectory is made of a single bang arc). 
 We want to show that
\bi
\iii[(a)]
When $\ph_{1}(0) \ph_{2}(0)> 0$  then the two switching functions are affine in a right-neighborhood of $0$.
\iii[(b)]
When $\ph_{1}(0) \ph_{2}(0)< 0$ the two switching functions are quadratic in a right-neighborhood of $0$.
\ei
On a bang arc the controls satisfy $|u_{1}|=|u_{2}|=1$ and thus  we can differentiate the identity \eqref{eq:pmpphimart0} and get
\begin{gather*}
\ddot{\ph}_{1}=-u_{2}\dot{\ph}_{3}=-4u_{2}(u_{1}-u_{2})\ph_{4},\quad
\ddot{\ph}_{2}=u_{1}\dot{\ph}_{3}=4u_{1}(u_{1}-u_{2})\ph_{4}.
\end{gather*}
In  case (a) we have that $u_{1}=u_{2}=\pm 1$, which implies $ \ddot{\ph}_{1}=\ddot{\ph}_{2}=0$.
In case (b) we have  $u_{1}-u_{2}=\pm 2$ and consequently  $\ddot{\ph}_{1}$ and $\ddot{\ph}_{2}$ are constant and nonzero (recall that $\ph_{4}\ne 0$ is constant).
The equations for case (a) are 
\begin{gather*}
\ph_{1}(t)=\ph_{1}(0)+t \dot{\ph}_{1}(0)=\ph_{1}(0)-u_{2}\ph_{3}(0)t,\\
\ph_{2}(t)=\ph_{2}(0)+t \dot{\ph}_{2}(0)=\ph_{2}(0)+u_{1}\ph_{3}(0)t,\\
\ph_{3}(t)=\ph_{3}(0).
\end{gather*}
Notice that $\ph_{3}(0)= y(0)\ph_{4}$ is zero if we start on the abnormal set. 
The equations for case (b) are
\begin{gather*}
\ph_{1}(t)=\ph_{1}(0)-u_{2}\ph_{3}(0)t-u_{2}(u_{1}-u_{2})\ph_{4}\frac{t^{2}}{2},\\
\ph_{2}(t)=\ph_{2}(0)+u_{1}\ph_{3}(0)t+u_{1}(u_{1}-u_{2})\ph_{4}\frac{t^{2}}{2},\\
\ph_{3}(t)=\ph_{3}(0)+(u_{1}-u_{2})\ph_{4}t.
  \end{gather*}
In particular, the constant $\ph_{4}$ determines the convexity of the quadratic arc of the switching functions.
\bl \label{l:singmart0} A regular bang arc can enter in a singular arc  only if the switching function is quadratic and has vanishing derivative at the switching point.    
\el
\begin{proof} Assume, for instance, that at some time $t_{0}\in I$ we have $\ph_{1}(t_{0})=1$ and $\ph_{2}(t_{0})=0$. Then the control $u_{1}(t)=\text{sign} \ph_{1}(t)$ is constantly equal to 1 in a neighborhood $U_{t_{0}}$ of $t_{0}$ and since $\ph_{3}$ is continuous we deduce that $\dot \ph_{2}=u_{1}\ph_{3}$ is also continuous in $U_{t_{0}}$. Since on the singular arc $\dot \ph_{2}=0$, we conclude.
\end{proof}
Next we discuss the possible behavior of the switching functions for regular arcs.
Let us 
assume that $\ph_{1}(0)>0$ and $\ph_{2}(0)<0$. In particular $\ph_{1}$ and $\ph_{2}$ are quadratic on a right-neighborhood of $0$.

\begin{figure}[h]
\begin{center}
	\scalebox{0.7}{\input{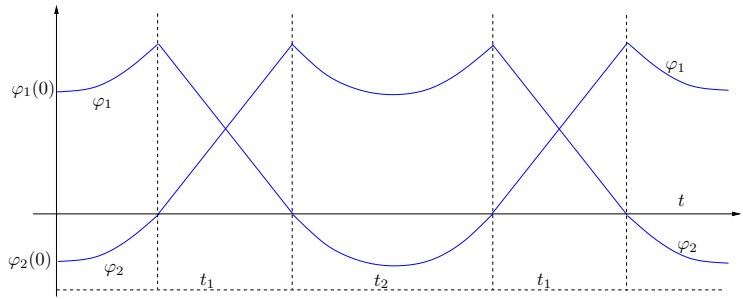}}
	\caption{Switching functions for the Martinet structure \eqref{eq:mart0} when $\ph_1$ is of type NI.}
\label{figura9}
\end{center}
\end{figure}

\begin{figure}[h]
\begin{center}
	\scalebox{0.75}{\input{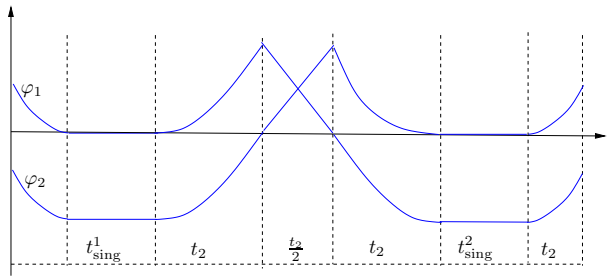}}
	\caption{Switching functions for the Martinet structure \eqref{eq:mart0} when $\ph_1$ is of type T. 
	The relation between the length of the third and fourth bang arcs can be easily deduced from the expression of the switching functions.}
\label{figura10}
\end{center}
\end{figure}

\begin{figure}[h]
\begin{center}
	\scalebox{0.7}{\input{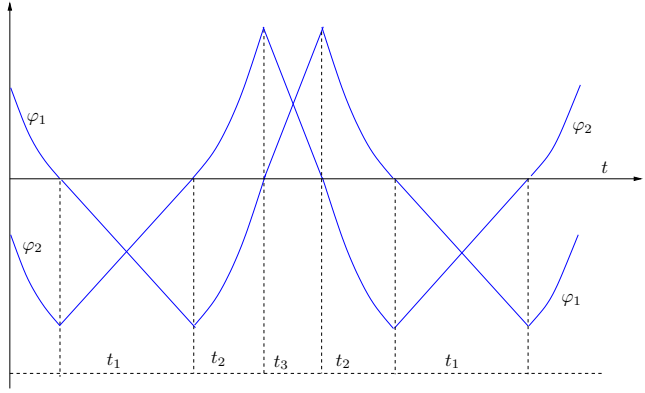}}
	\caption{Switching functions for the Martinet structure \eqref{eq:mart0} when $\ph_1$ is of type I.}
\label{figura11}
\end{center}
\end{figure}

We are reduced to three possible cases for the the switching function $\ph_{1}$: 
\bi
\iii[-] it never vanishes in the quadratic part (we say that $\ph_1$ is of \emph{type NI}, for \emph{not intersecting}), 
\iii[-] it vanishes 
in the quadratic part and is tangent to the zero level (\emph{type T} for \emph{tangent}), 
\iii[-] it vanishes 
in the quadratic part and is transversal to the zero level (\emph{type I} for \emph{intersecting}).
\ei
In Figure \ref{figura9} we picture the switching functions when $\ph_1$ is of type NI, while Figures~\ref{figura10} and  \ref{figura11} correspond to  type T and  type I, respectively. 


Assuming that there are only regular bang arcs along the trajectory (as it is always the case when $\ph_1$ is of type NI or I) we have the following result.
\bp \label{p:onlyr} The switching functions of a trajectory that has only regular bang arcs are periodic. 
\ep
The proof of Proposition \ref{p:onlyr} is a simple consequence of the formulas of the switching functions and Lemma \ref{l:singmart0}. 
When $\ph_1$ is  of type T, the only freedom is in the length of singular arcs.
The order in which the switching occur is as in Figures~\ref{figura9}, \ref{figura10} and  \ref{figura11}, up to the symmetry 
which sends $(\ph_1,\ph_2,\ph_3,\ph_4)$ into $(-\ph_1,-\ph_2,\ph_3,-\ph_4)$ (which corresponds to a  reflection $y\to -y$). 
\brem
It is easy to see from equations \eqref{eq:pmpphimart0} and Remark~\ref{r:ph4} that if  $\ph_1$ is of type I 
then the $y$ coordinate of the corresponding trajectory has constant sign. When restricting our attention to trajectories starting on the plane $\{y=0\}$, we can then
exclude that $\ph_1$ is of type I.  
\erem


\subsubsection{Bound on the number of  regular arcs for optimal trajectories}

The goal of this section is to prove the following result.

\bp \label{lemma-bang-optimal-martinet1}
A  bang-bang trajectory with at least one regular arc and with more than $7$  arcs (either bang or singular) is not optimal.
If, moreover, the trajectory starts on the plane $\{y=0\}$ and has more than $5$ arcs, then it is not optimal. 
\ep

The proof works by applying several times Theorem~\ref{thm2nd}. We should distinguish trajectories for which the switching functions are of one of the three types NI, T, and I. 

In order to reduce the number of cases to be studied, we use the fact that  time-reversion and reflection  $y\to -y$ lead to trajectories with equivalent optimality properties.

\subsubsection*{Switching functions of type NI}

We start by considering $\ph_1$ of the type NI, as in Figure~\ref{figura9}.

\bl \label{lemma-bang-optimal-martinet1-NI}
A regular bang-bang trajectory of type NI with more than $5$  arcs is not optimal.
If, moreover, the trajectory starts on the plane $\{y=0\}$ and has more than $3$ arcs, then it is not optimal. 
\el
\begin{proof}
We prove the first part of the lemma
by showing that concatenations of the type 
\begin{equation}\label{conca-NI}
 (1,-1),(1,1),(-1,1),(1,1),(1,-1)
 \end{equation}
are not optimal. All concatenations of $6$ bang arcs, indeed, contain a concatenation of this type, up to symmetries (see Figure~\ref{figura12}).

For concatenations of type \eqref{conca-NI}, applying Theorem~\ref{thm2nd} at the second switching time, we get by computations as the one seen in the previous sections that the space $W$ and 
the quadratic form $Q$ 
in the statement of  Theorem~\ref{thm2nd} can be written as
$$W=\{(\al_0,\al_1,0,-\al_1,-\al_0)\mid \al_0,\al_1\in\R\},\quad Q(\al_0,\al_1)=8(t_1\al_0^2+t_2 \al_0\al_1).$$
Since $Q$  is  not negative semidefinite, the corresponding trajectory is not optimal.

In order to conclude the proof of Lemma~\ref{lemma-bang-optimal-martinet1-NI}, notice that, by metric considerations, if the trajectory starts from the plane $\{y=0\}$, then it stops to be optimal at the middle of the third bang arc, see Figure \ref{figura12}. 
\end{proof}
%
%

%

\subsubsection*{Switching functions of type T}

We prove here the following result concerning
 trajectories corresponding to switching functions of the type T as in Figure~\ref{figura10}.

\bl \label{lemma-bang-optimal-martinet1-T}
A trajectory of type T with more than $7$  arcs is not optimal.
If, moreover, the trajectory starts on the plane $\{y=0\}$ and has more than $5$ arcs, 
then it is not optimal. 
\el
\begin{proof}
We first  consider the situation where $t^k_{\mathrm{sing}}>0$ for every $k$.
We notice that 
every concatenations of $8$ arcs contains, up to symmetries, a concatenation of $6$ arcs of the type
\begin{equation}\label{conca-T-2}
(1,-1),(1,1),(-1,1),(-1,-1),(-1,1),(1,1).
\end{equation}
(See Figure~\ref{figura13}).
We are going to  show that a concatenation as in \eqref{conca-T-2} is not optimal.

For concatenations of type \eqref{conca-T-2}, applying Theorem~\ref{thm2nd} at the third switching time (at which $y=0$), we get that
the space $W$ and 
the quadratic form $Q$ 
in the statement of  Theorem~\ref{thm2nd} are written as
\begin{gather*}
W=\{(\al_0,\al_1,\al_2,-\al_0,\al_0-\al_2,-\al_0-\al_1)\mid \al_0,\al_1,\al_2\in\R\},\\
Q(\al_0,\al_1,\al_2)=
2( t_2 - 2 t^2_{\mathrm{sing}}) \al_0^2+8t_2 \al_0\al_1+8t^2_{\mathrm{sing}} \al_0\al_2 - 4 t^2_{\mathrm{sing}} \al_2^2. 
\end{gather*}
Notice that $Q$  is  not negative semidefinite, since $Q(\eps,1/\eps,0)=2\eps^2 ( t_2 - 2 t^2_{\mathrm{sing}})+8t_2 >0$ for $\eps$ small enough. Hence, the corresponding trajectory is not optimal. 

In the case where $t^2_{\mathrm{sing}}=0$,
a concatenation as in  \eqref{conca-T-2} reduces to  a concatenation of $4$ bang arcs
$$ (1,-1),(1,1),(-1,1),(1,1).$$
Considering the following arc, we recover a concatenation 
as in \eqref{conca-NI}, for which the same computations as in the previous section show non-optimality. 

The proof of Lemma~\ref{lemma-bang-optimal-martinet1-T} can be concluded as before by metric considerations for trajectories starting from the plane $\{y=0\}$, see Figure \ref{figura13}. 
%
 \end{proof}

\subsubsection*{Switching functions of type I}

We consider here trajectories corresponding to switching functions of type I as in Figure~\ref{figura11}.
Notice that such trajectories never cross the plane $\{y=0\}$.

We  prove the following result.
 
\bl \label{lemma-bang-optimal-martinet1-I}
A regular bang-bang trajectory of type I with more than $5$  arcs is not optimal.
\el
\begin{proof}

By the same symmetry considerations as in the cases NI and T, we are left to prove that 
 concatenations of the type 
\begin{equation}\label{conca-I-1}
 (1,-1),(-1,-1),(-1,1),(1,1),(1,-1),(-1,-1)
 \end{equation}
and
\begin{equation}\label{conca-I-2}
(-1,1),(1,1),(1,-1),(-1,-1), (-1,1),(1,1)
 \end{equation}
are not optimal (see Figure~\ref{figura14}). Notice than in both cases the 
 trajectory is  contained in $\{y<0\}$. 

The application of Theorem~\ref{thm2nd}  to the two cases is very similar leading (computing the quadratic form $Q$ at the second switching time $\tau_2$) to the expressions
$$Q(\al_0,\al_1,\al_2)=-4 (t_1 - t_3) \al_0^2   
 -4 y(\tau_2) \al_0\al_1 - 4 (-2 t_2 + 2 t_3 +  y(\tau_2)) \al_0\al_2- 4 (2 t_2 - t_3 - y(\tau_2)) \al_2^2$$
and
$$Q(\al_0,\al_1,\al_2)=-4 (t_1 - t_3) \al_0^2   
 -4 y(\tau_2) \al_0\al_1 - 4 (-2 t_1 + 2 t_2 +  y(\tau_2)) \al_0\al_2- 4 ( t_1 - 2t_2 - y(\tau_2)) \al_2^2$$
respectively.
In both cases, since $y(\tau_2)< 0$, one has that $Q(\eps,1/\eps,0)=-4 y(\tau_2)+O(\eps^2)$ is positive for $\eps$ small enough.
Theorem~\ref{thm2nd}  then allows to conclude that 
 the corresponding trajectories are not optimal. 
\end{proof}

 \subsubsection{Optimal trajectories and shape of the unit ball}
Here we present the different pictures for the $(x,y)$-components of trajectories corresponding to switching functions of the form NI, T and I. The dashed lines correspond to the part of the trajectory which is no more optimal. 

\begin{figure}[h]
\begin{center}
	\scalebox{.7}{\input{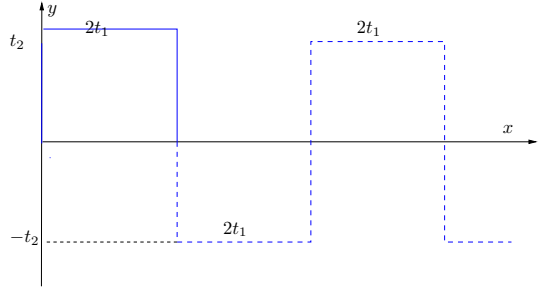}}
	\caption{Regular bang-bang  trajectories of type NI for the Martinet structure \eqref{eq:mart0}.}
\label{figura12}
\end{center}
\end{figure}

\begin{figure}[h]
\begin{center}
	\scalebox{.75}{\input{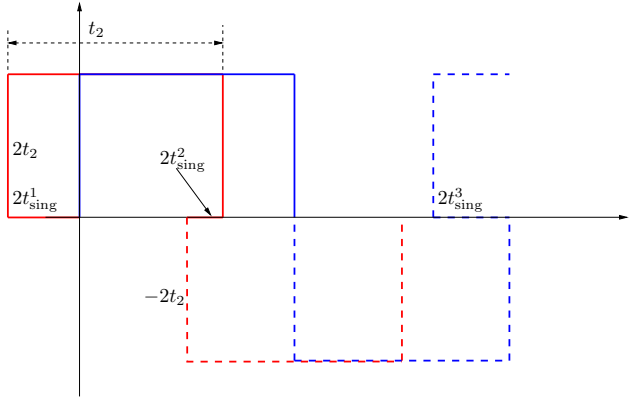}}
	\caption{Regular/singular bang trajectories of type T for the Martinet structure \eqref{eq:mart0}. }
	  \label{figura13}
\end{center}
\end{figure}

\begin{figure}[h]
\begin{center}
	\scalebox{.7}{\input{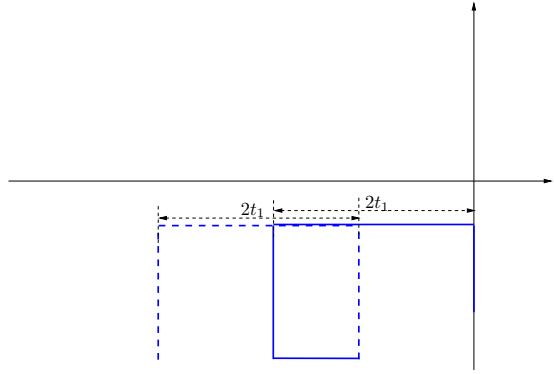}}
	\caption{Regular bang-bang trajectories of type I for the Martinet structure \eqref{eq:mart0}.}
\label{figura14}
\end{center}
\end{figure}

\begin{figure}[h]
\begin{center}
\scalebox{.25}{\includegraphics{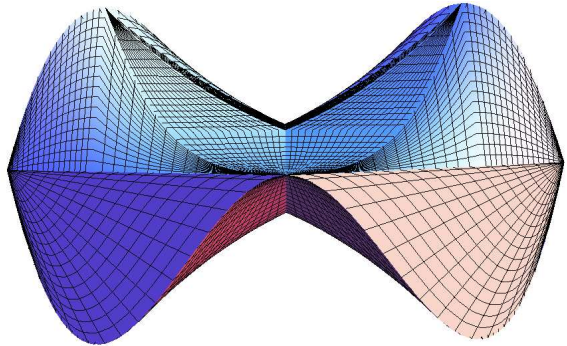}}
	\caption{Unit sphere for the Martinet structure \eqref{eq:mart0}, view from the $x$-axis.}
\label{figura15A}
\end{center}
\end{figure}

\begin{figure}[h]
\begin{center}
\scalebox{.25}{\includegraphics{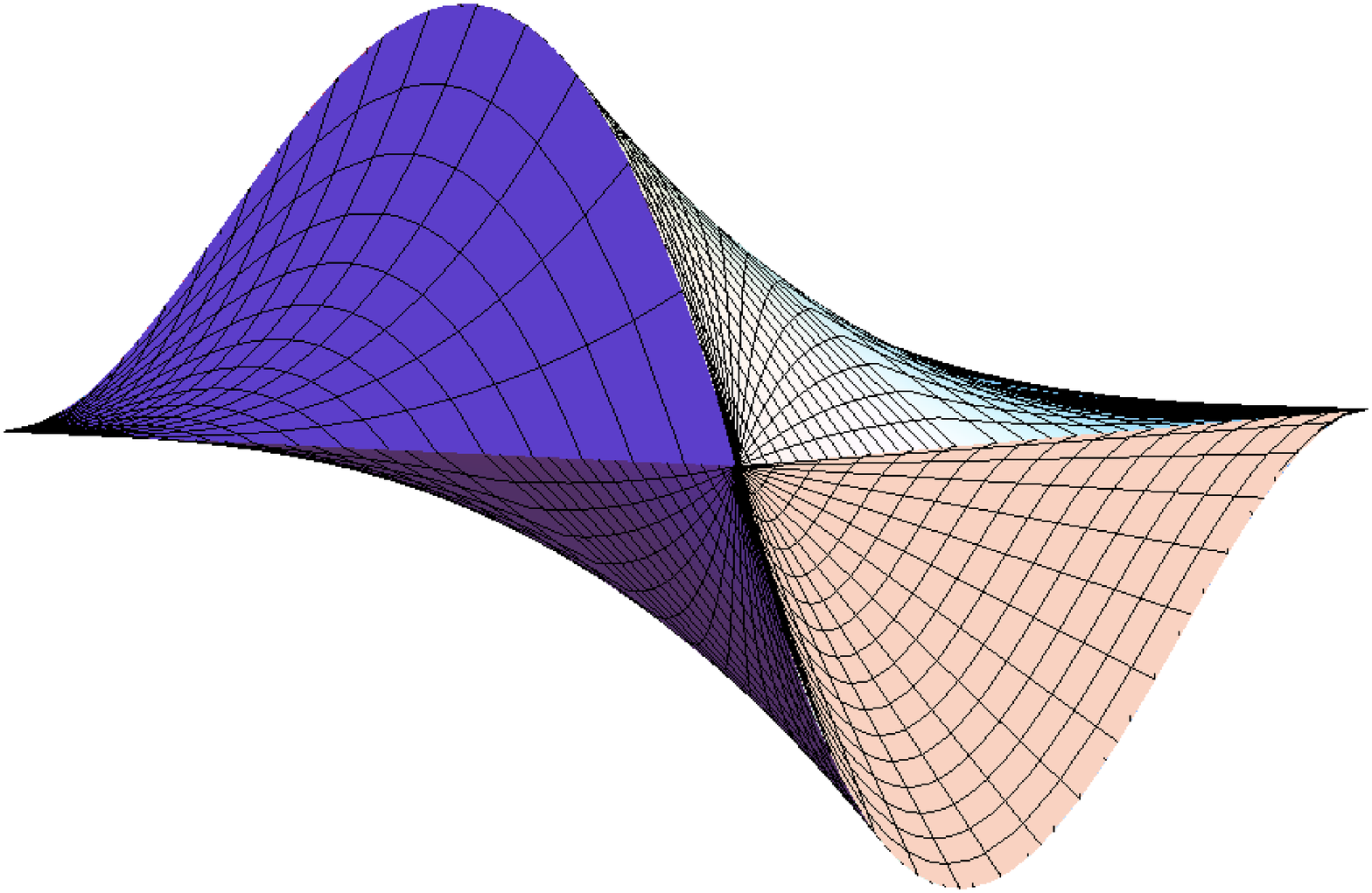}}
	\caption{Unit sphere for the Martinet structure \eqref{eq:mart0}, view from the $y$-axis.}
\label{figura15B}
\end{center}
\end{figure}

In the case of trajectories of type NI we have the behavior in Figure \ref{figura12}.
Trajectories of type T have singular arcs of arbitrary length (see two examples  in Figure \ref{figura13}). Notice that the switching to singular always happens at points where $y=0$, namely on the Martinet surface. 
The last case is given by trajectories of type I (see  Figure \ref{figura14}). In this case the trajectory is contained in a strip $y_{0}\leq y(t)\leq y_{1}$ with eiter $0<y_{0}$ or $y_{1}>0$. 
In view of the optimality results one gets the following picture of the unit ball in the Martinet structure \eqref{eq:mart0}, see Figures \ref{figura15A} and \ref{figura15B}.

\subsection{The second structure}
The second  sub-Finsler Martinet structure on  $\mb{R}^{3}$ that we are going to consider is the sub-$\ell^\infty$ structure
determined by the vector fields 
\bqn\label{eq:mart00}
X_{1}=\partial_{x}+y^{2} \partial_{z},\qquad X_{2}=\partial_{y}.
\eqn
We introduce the  vector fields 
\begin{equation}\label{BrRe2}
X_{3}:=[X_{1},X_{2}]=2y\partial_{z},\quad X_{4}:=[X_{1},[X_{1},X_{2}]]=0,\quad X_{5}:=[X_{2},[X_{1},X_{2}]]=2\partial_{z},
\end{equation}
and the switching functions 
$$\varphi_{i}(t)=\la \lambda(t) ,X_{i}(\gamma(t))\ra,\qquad i=1,\ldots,5.$$
The functions $\ph_{i}$ satisfy the following system of differential equations
\begin{gather} \label{eq:pmpphimart00}
\dot{\ph}_{1}=-u_{2}\ph_{3},\qquad \dot{\ph}_{2}=u_{1}\ph_{3}, \qquad \dot{\ph}_{3}=u_{1}\ph_{4}+u_{2}\ph_{5},
\qquad \dot{\ph}_{4}=0,\qquad \dot{\ph}_{5}=0.
\end{gather}
In this case the additional relations given by the bracket relations \eqref{BrRe2} are  $\ph_{3}=y\ph_{5}$ 
and $\ph_{4}=0$. In particular system \eqref{eq:pmpphimart00} reduces to
\bqn \label{eq:pmpphimart000}
\dot{\ph}_{1}=-u_{2}\ph_{3},\qquad \dot{\ph}_{2}=u_{1}\ph_{3}, \qquad \dot{\ph}_{3}=u_{2}\ph_{5},\qquad 
 \dot{\ph}_{5}=0.
\eqn
Reasoning as in Lemma \ref{l:reasoning}, we have the following characterization of abnormal arcs. 
\bl 
The nontrivial abnormal arcs on the sub-$\ell^{\infty}$ structure on $\R^{3}$ defined by the vector fields \eqref{eq:mart00} are the horizontal lines contained in the plane $\{y=0\}$.
\el
Indeed  abnormal trajectories are described by the equations $\{y=0,z=z_{0}\}$. The fact that these trajectories are the same in the two Martinet structure under consideration reflects the fact that abnormal trajectories are independent of the choice of the frame (they depend only on the distribution).
%

\subsubsection{Singular arcs}
The situation is not in this case symmetric with respect to $\ph_{1}$ and $\ph_{2}$ singular.
Let us first consider a $\ph_{1}$-singular arc, i.e., $\ph_{1}\equiv0$ on $I$. Since $u_{2}=\pm 1$ is constant and $0\equiv \dot \ph_{1}=-u_{2}\ph_{3}$, it follows that  $\ph_{3}(t)= 0$ for all $t\in I$. From $\dot\ph_{3}=u_{2}\ph_{5}$ we deduce that $\ph_{5}=0$ and $u_{1}$ is arbitrary. In particular $\ph_3$ is identically equal to zero, which implies that the trajectory stays singular for all times.

Let us then consider a $\ph_{2}$-singular arc, i.e., $\ph_{2}\equiv0$ on $I$. Since $u_{1}=\pm 1$ is constant and $0\equiv \dot \ph_{2}=u_{1}\ph_{3}$ it follows that  $\ph_{3}(t)= 0$ for all $t\in I$. We deduce that $0=\dot\ph_{3}=u_{2}\ph_{5}$. Hence, either $\ph_{5}=0$ with $u_{2}$ arbitrary, or $u_{2}=0$ on $I$ with $\ph_{5}\neq 0$. In the first case the trajectory stays singular for all times, in the second case the singular arc is contained in the plane $\{y=0\}$ and coincides with an abnormal trajectory on the interval $I$.

\subsubsection{Regular arcs}
The analysis is similar to that of the previous section. We therefore omit the computations, which yield the following result.

\bp \label{lemma-bang-optimal-martinet2}
A  regular trajectory with more than 6  arcs (either bang or singular) is not optimal.
\ep

\begin{figure}[h]
\begin{center}
	\includegraphics[scale=0.3]{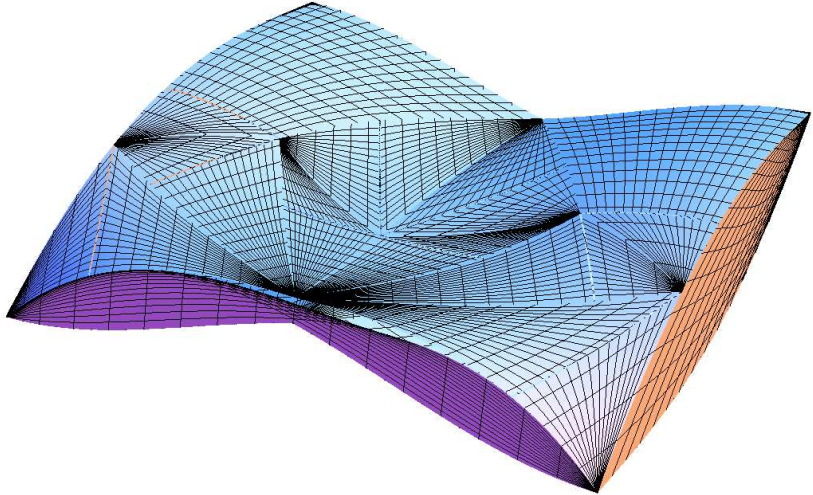}
	\caption{Unit sphere for the Martinet structure \eqref{eq:mart00}.}
\label{figura20}
\end{center}
\end{figure}

 The picture of the unit ball in the Martinet structure \eqref{eq:mart00} is given in Figure \ref{figura20}.

\section{Euclidean rectifiability and semi-analyticity of  spheres}

By construction the $\ell^\infty$ spheres for  the Heisenberg, Grushin and Martinet structures studied above are homeomorphic to Euclidean spheres ($S^2$ for Heisenberg and Martinet  and $S^1$ for Grushin). Moreover these spheres are graphs of  piecewise-polynomial functions. 

It follows that these spheres are Euclidean rectifiable and semi-analytic. Theorem \ref{222} follows.

{\small
\bibliography{subfinsler_biblio_DONOTMODIFY}
\bibliographystyle{amsalpha}

}

\end{document}